\documentclass[11pt]{article}
\usepackage{amsmath}
\usepackage{amssymb}
\usepackage{amsfonts}
\usepackage{mathrsfs}
\usepackage{shortvrb,psfrag}
\usepackage{enumerate}
\usepackage{color}
\usepackage[dvips]{graphicx}

\let\al=\alpha
\let\be=\beta

\let\f=\frac
\let\p=\psi
\let\veps=\varepsilon

\let\Lam=\Lambda

\let\Om=\Omega

\def\cC{{\cal C}}

\def\cE{{\cal E}}
\def\cF{{\cal F}}
\def\cG{{\cal G}}

\def\no{\noindent}
\def\na{\nabla}
\def\p{\partial}
\def\eqdef{\buildrel\hbox{\footnotesize def}\over =}
\def\endproof{\hphantom{MM}\hfill\llap{$\square$}\goodbreak}

\newcommand{\beq}{\begin{equation}}
\newcommand{\eeq}{\end{equation}}
\newcommand{\ben}{\begin{eqnarray}}
\newcommand{\een}{\end{eqnarray}}
\newcommand{\beno}{\begin{eqnarray*}}
\newcommand{\eeno}{\end{eqnarray*}}

\newcommand{\vcx}{\vec{x}}
\newcommand{\vcu}{\vec{u}}
\newcommand{\ud}{\mathrm{d}}

\newcommand{\uu}{\mathbf{u}}
\newcommand{\vv}{\mathbf{v}}
\newcommand{\ww}{\mathbf{w}}

\newcommand{\nn}{\mathbf{n}}
\newcommand{\nna}{\mathbf{n}_{\alpha}}
\newcommand{\nnb}{\mathbf{n}_{\beta}}

\newcommand{\RR}{\mathbf{R}}
\newcommand{\RT}{\mathbf{T}}
\newcommand{\II}{\mathbf{I}}
\newcommand{\RRa}{\mathbf{R}_{\alpha}}
\newcommand{\RRb}{\mathbf{R}_{\beta}}
\newcommand{\TR}{\widetilde{\mathbf{R}}}
\newcommand{\HR}{\widehat{\mathbf{R}}}
\newcommand{\TRa}{\widetilde{\mathbf{R}}_{\alpha}}
\newcommand{\HRa}{\widehat{\mathbf{R}}_{\alpha}}
\newcommand{\TRb}{\widetilde{\mathbf{R}}_{\beta}}
\newcommand{\HRb}{\widehat{\mathbf{R}}_{\beta}}

\newcommand{\FF}{\mathfrak{F}}
\newcommand{\ta}{\mathbf{t}^1}
\newcommand{\taa}{\mathbf{t}^1_{\alpha}}
\newcommand{\tab}{\mathbf{t}^1_{\beta}}
\newcommand{\tb}{\mathbf{t}^2}
\newcommand{\tba}{\mathbf{t}^2_{\alpha}}
\newcommand{\tbb}{\mathbf{t}^2_{\beta}}
\newcommand{\aaa}{\mathbf{a}}

\newcommand{\aab}{\mathbf{a}_{\beta}}
\newcommand{\se}{\sqrt{E}}
\newcommand{\Ea}{E_{\alpha}}
\newcommand{\Eb}{E_{\beta}}
\newcommand{\un}{U^{n}}
\newcommand{\una}{U^{n}_{\alpha}}
\newcommand{\unb}{U^{n}_{\beta}}
\newcommand{\uua}{U_1}
\newcommand{\uub}{U_2}
\newcommand{\uupa}{\uu_{\alpha}}
\newcommand{\uupb}{\uu_{\beta}}
\newcommand{\uaa}{U_{1\alpha}}
\newcommand{\uab}{U_{1\beta}}
\newcommand{\uba}{U_{2\alpha}}
\newcommand{\ubb}{U_{2\beta}}
\newcommand{\Wa}{W_1}
\newcommand{\Wb}{W_2}

\newcommand{\ab}{\al\be}

\newtheorem{Theorem}{Theorem}[section]

\newtheorem{Proposition}[Theorem]{Proposition}
\newtheorem{Lemma}[Theorem]{Lemma}

\newtheorem{Remark}[Theorem]{Remark}

\begin{document}
\title{Well-posedness of Hydrodynamics on the Moving Elastic Surface}

\author{Wei Wang,\, Pingwen Zhang \,and \,Zhifei Zhang\\[2mm]
{\small LMAM and School of  Mathematical Sciences, Peking University, Beijing 100871, China}\\
{\small E-mail: wangw07@pku.edu.cn, pzhang@pku.edu.cn,
zfzhang@math.pku.edu.cn}}

\date{\today}
\maketitle

\begin{abstract}
The dynamics of a membrane is a coupled system comprising a moving elastic
surface and an incompressible membrane fluid. We will consider a reduced elastic surface
model, which involves the evolution equations of the moving surface, the dynamic
equations of the two-dimensional fluid, and the incompressible equation,
all of which operate within a curved geometry. In this paper, we prove
the local existence and uniqueness of the solution to the
reduced elastic surface model by reformulating the model into a new system in the
isothermal coordinates. One major difficulty is that of constructing an appropriate iterative scheme
such that the limit system is consistent with the original system.

\end{abstract}

\setcounter{equation}{0}
\section{Introduction}
This paper is concerned with the hydrodynamics on the moving
surface of bio-membrane, which as the outerwear of living cells and organelles
plays an important role in the life process. Consisting of lipids,
proteins and carbohydrates, the structures and
properties of bio-membrane are very complex. In general, bio-membrane can be viewed
as a 2-dimensional fluid surface consisting of a lipid bilayer, as
the lipid molecules can move freely on the surface but cannot
escape from it. The fluid is viscous and can be viewed as incompressible
because it typically has a large tensile module.
Moreover, this 2-dimensional fluid surface is bend-resistent. Hence, it tends to
minimize the Helfrich energy under the fixed area condition (guaranteed by the incompressible condition)
\begin{equation}\label{helfrich-e}
E_H=\int\big(c_1(H-B)^2-c_2K\big)\ud{\sigma},
\end{equation}
where $H$ and $K$ are the mean curvature and the Gaussian curvature, respectively, $B$ is the spontaneous curvature
that reflects the initial or intrinsic curvature of the membrane, $c_1$ and $c_2$ are the elastic coefficients,
and $\ud\sigma$ is the area form of the surface \cite{DHelfrich}. When $c_2$ is uniform on the membrane,
$\int K\ud{\sigma}$ is a constant determined by the topology of the membrane.
When $B \equiv0$, $E_H$ is called the Willmore energy in geometry.
A number of studies based on Helfrich's bending energy model explore the mechanics
of bio-membrane, for example, see \cite{Steig0, CapGuv, LM}.

During the past several decades, membrane dynamics  have
received considerable attention. Researchers from different fields have developed several models
with/without the surrounding fluid to study the behaviors of
the membrane.  For the models without surrounding fluid, see \cite{ Waxman, Waxman2, Scriven, Steig,
CapGuv}, and for the models with surrounding fluid, see \cite{Poz, Poz6, Miao, Pow}.

Waxman \cite{Waxman} may have been the first to study the dynamics of bend-resistant
bio-membrane using a model without surrounding fluid and in which the incompressibility, bend-resistance, and
viscosity effects are all considered.  However, Waxman's model does not preserve the energy
dissipation law. In \cite{Hu-Zhang}, Hu-Zhang-E introduced a director field to
represent the direction of lipid molecules at every material point
and developed an elastic energy model based on the Frank energy of
the smectic liquid crystal. When the director is constrained to
the normal of the surface, they obtain a reduced elastic surface
model, that is very close to Waxman's model, but adds one term to the
in-plane stresses whereby the model satisfies a natural energy
dissipation law. In the elastic surface model, the dynamics of the membrane
involves the evolution equations of the moving surface, the dynamic
equations of the two-dimensional fluid, and the incompressible equation,
all of which operate within a curved geometry.

For a surface membrane $\Gamma=\RR(u_1,u_2,t)$,
we denote by $\aaa_{\alpha}$ the tangent vectors of $\Gamma$, $\nn$  the
unit normal vector, $(a_{\alpha\beta})_{1\le \al,\be\le 2}$ the covariant metric tensor,
$\Delta_\Gamma$ the Lapalace-Beltrami operator, $K$  the
Gaussian curvature, and $H$  the mean curvature. In the simple case, the reduced elastic surface
model takes the following form:
\begin{eqnarray}\label{eq:reduced model}
\left\{
\begin{array}{l}
\f {\p\RR} {\p t}=\vv(u_1,u_2,t),\\
\frac{\partial{\vv}}{\partial{t}}=(-\Pi{a^{\alpha\beta}\aaa_{\alpha}})_{,\beta}+2\varepsilon_0
(S^{\al\be}\aaa_{\alpha})_{,\beta}
-\frac{1}{2}\big(\Delta_{\Gamma}{H}+2H(H^2-K)\big)\nn,\\
\na_\Gamma\cdot \vv=0.
\end{array}\right.
\end{eqnarray}
Here, $\vv$ is the velocity of the fluid,  $\Pi$ is the surface
pressure, $S^{\al\be}$ is the rate of the surface strain, and the constant $\varepsilon_0>0$ is the shear viscosity.
The notation $()_{,\beta}$ denotes the covariant derivative.
The first term on the right-hand side of the second equation is induced by the incompressible condition $\na_\Gamma\cdot \vv=0$, and the
surface pressure $\Pi$ can be viewed as a Lagrangian multiplier;
the second term describes the viscosity of the fluid on the surface;
the third term is the elastic stress induced by the Helfrich bending energy
(\ref{helfrich-e}) with $B=0$. Please see Section 2 or \cite{Hu-Zhang}
for more detail.

When the interaction with bulk fluid is considered, Hu-Zhang-E
\cite{Hu-Zhang} also derived the incompressible membrane-fluid coupling
system in the form
\begin{eqnarray}\label{coupling1-1}
\left\{
\begin{array}{rcll}
\uu_t+\uu\cdot\nabla\uu&=&-\nabla p+\nu\Delta\uu, \quad &{\textrm{in}} ~\Omega,\\
\nabla\cdot\uu&=&0,\quad &\textrm{in} ~\Omega,\\
~[-p\II+\tau]\cdot\nn&=&\mathbf{F},\quad&\textrm{on}~\Gamma,\\
~[\uu]&=&0,\quad&\textrm{on}~\Gamma,\\
~\nabla_{\Gamma}\cdot\uu&=&0,\quad\quad\quad\quad\quad\quad&\textrm{on}~\Gamma,
\end{array}\right.
\end{eqnarray}
where $\tau=\nu(\nabla\uu+\nabla\uu^T)$ is the stress of the bulk
fluid, $\mathbf{F}$ is given by the right-hand side of the second equation of
(\ref{eq:reduced model}), $\Om$ is the fluid domain, $\Gamma$ is
the time-dependent surface of the membrane included in $\Om$,
and $[\cdot]$ denotes the jump across the membrane. In a recent review
paper \cite{Pow}, a similar model was derived via the direct variational
method. Compared with the classical free boundary problem of the Navier-Stokes equations,
the main difference is that the system (\ref{coupling1-1})
contains two unknown pressures: the pressure $p$ of the
surrounding fluid and the pressure $\Pi$ of the
membrane defined on the surface, where $\Pi$ is determined by
the incompressible condition $\nabla_{\Gamma}\cdot\uu=0$.
Due to the coupling between $p$ and $\Pi$, solving the membrane-fluid coupling system
(\ref{coupling1-1}) is still challenge, both mathematically and
numerically. In some specific case (e.g., when the velocity of the surrounding
fluid is small), the main influence of the bulk fluid is to maintain
the enclosed volume of the membrane. In such cases for simplicity, it can be
replaced by introducing osmotic pressure. Moreover, although the reduced model (\ref{eq:reduced model}) neglects the fluid interaction,
numerical simulation \cite{HuPhd} also convinces us that this model can be used
to reconstruct some important physical processes, such as
exocytosis and endocytosis.

To our knowledge, few mathematical results such as the
well-posedness for the fluid bio-membrane dynamics are available. In
\cite{CCShkoller}, Cheng-Coutand-Shkoller studied the bulk fluid
interacting with a membrane considered a nonlinear elastic
bio-fluid shell and modeled by the nonlinear Saint Venant-Kirchhoff
constitutive law, where the membrane is compressible and the surface
fluid is inviscid. In \cite{HuSongZhang}, Hu-Song-Zhang
proved the local existence and uniqueness of (\ref{eq:reduced model}) for a simplified case when
the membrane is cylindrical. In this case,  the membrane is
similar to a 1-D incompressible string such that the fluid vanishes. With the introduction of  the arc length parameter
and the tangent angle of the curve,
the system is transformed into a fourth-order wave equation for the tangent angle $\al$
coupled with an elliptic equation:
\beno
\left\{\begin{array}{l}
\al_{tt}=g_1+2T_s\al_s+T\al_{ss}-(B+\al_s)_{sss}+\al_s^2(B+\al_s)_s,\\
-T_{ss}+\al_s^2T=g_2+\al_t^2+2(B+\al_s)_{ss}\al_s+(B+\al_s)_s\al_{ss},
\end{array}\right.
\eeno
where $g_1, g_2$, and $B$ are the given smooth functions.

The purpose of this paper is to prove the local well-posedness of
the system (\ref{eq:reduced model}). This is also a key step toward
understanding and solving the membrane-fluid coupling system
(\ref{coupling1-1}). Our result is stated as follows.

\begin{Theorem}\label{thm:main}
Let $s=2k$ for some integer $k\ge 3$.
Assume that the initial velocity $\vv_0\in H^{s-1}$ and the initial closed surface $\RR_0\in H^{s+1}$.
There exists $T>0$ such that the system (\ref{eq:reduced model})
has a unique solution $(\vv(t), \RR(t))$  on $[0,T]$ satisfying
\beno
\vv\in C([0,T];H^{s-1}),\qquad \RR\in C([0,T];H^{s+1}).
\eeno
\end{Theorem}

\begin{Remark}
The regularity we imposed on the initial data should not be
optimal. To simplify the analysis, we will work in a functional space with
high regularity.
\end{Remark}

System (\ref{eq:reduced model}) is a coupled system of
parabolic, hyperbolic, and elliptic equations. The evolution
equations of the tangential velocities are parabolic, the evolution
equations of the normal velocity and the mean curvature constitute a
hyperbolic system, and the pressure satisfies an elliptic equation,
see (\ref{eq:full system-R})-(\ref{eq:full system-velocity free}).
Because the surface is moving, it seems natural to solve (\ref{eq:reduced
model}) in the framework of Lagrangian coordinates. However, some essential difficulties
will arise. Let us explain it in what follows.

Assume that the initial velocity $\vv_0\in H^{s-1}$ and the initial surface $\RR_0\in H^{s+1}$.
Because the tangential velocity $v^\al$ satisfies the parabolic equation, and
the normal velocity $v^n$  and the mean curvature $H$ together satisfy the
hyperbolic system, it seems natural to expect
$v^\al$ to belong to $L^2(0,T;H^{s})$, and $(v^n, H)$ to belong to $L^\infty(0,T;H^{s-1})$.
However, these estimates depend on the $H^{s}$ regularity of the metric of the surface.
Hence, we have to recover the $H^{s}$ regularity of the metric from $(v^\al, v^n, H)$
in order to close the energy estimates.
In the Lagrangian coordinates, we have
\beno \RR_t=\vv(u_1,u_2, t),
\eeno
which tells us that $\RR \in L^\infty(0,T;H^{s-1})$ by the estimate for the velocity.
Hence the metric has only $H^{s-2}$ regularity (a loss of two derivatives). Maybe, one wants to use
the regularity of the mean curvature to gain the regularity of
$\RR$(Note that formally, $H^{s-1}$ regularity of the mean curvature suggests that the free surface has $H^{s+1}$
regularity). However, we cannot expect  $\RR$ to have
more regularity in the Lagrangian coordinates, see the example and argument of Section
5 in \cite{Shatah}.

Another way to solve the system is to
represent the moving surface locally by $x_3=g(x_1,x_2,t)$, where $g$
satisfies the following hyperbolic equation
\beno \f 1
{\sqrt{1+|\na_\Gamma g|^2}}g_{tt}+\Delta_\Gamma\Big(\textrm{div}_\Gamma(\f
{\na_\Gamma g} {\sqrt{1+|\na_\Gamma g|^2}})\Big) =\textrm{lower-order terms}.
\eeno
However, if we make an energy estimate for this equation,
the estimate is also not closed, since the lower-order
terms contain the third-order derivative of $g$, which cannot be
controlled by the main part.

Motivated by \cite{Ambrose}, we will use the
isothermal coordinates to re-parameterize the surface. There are two main
advantages adopting the isothermal coordinate: (1) we can gain two
more regularities for the surface from the regularity of the mean
curvature, and (2) the coefficients of the first fundamental form
have the same regularity as the surface.
Indeed, there are the following important relations between the surface $\Gamma=\RR(u_1,u_2)$,
the first fundamental form $E$, and the mean curvature $H$
when $(u_1,u_2)$ is taken as the isothermal coordinate of $\Gamma$:
\ben\label{eq:geomerty}
\Delta\RR=2EH\nn, \quad\Delta E=2\big(\p_{u_1}\p_{u_2}\RR\cdot\p_{u_1}\p_{u_2}\RR-\p_{u_1}^2\RR\cdot\p_{u_2}^2\RR\big).
\een
Here $\Delta=\p_{u_1}^2+\p_{u_2}^2$, and $\nn$ is the unit normal of $\Gamma$.

In general, it is difficult to construct an approximate system preserving the isothermal relation.
As the solution of the approximate system does not satisfy the important geometric relation (\ref{eq:geomerty}),
there will also be derivative loss once we make the energy estimates for the approximate system.
To overcome this difficulty, we incorporate the relation (\ref{eq:geomerty}) into our iterative scheme.
However, this produces another very troubling problem---one that arises for the most part from
the construction of the iterative scheme and relates to the equivalence of the two systems. The problem is
this: we do not know and need to establish whether  the limit system is equivalent
to the original system, and proving such equivalence involves very complicated geometric calculations.

This paper is organized as follows. In the next section, we review
some formulae for the evolving surfaces and introduce the reduced elastic surface model.
In Section 3, we derive an equivalent system in the isothermal
coordinate by decomposing the velocity into tangential
and normal components. Section 4 is devoted to studying the linearized system.
In Section 5, we prove our main results, including the
construction of the iteration scheme, nonlinear estimates, the
convergence of the iteration procedure, and the consistency between the limit system and the
original system.

\setcounter{equation}{0}
\section{The elastic model of an incompressible fluid membrane}
In this section, we provide a short derivation of the dynamic
model of an incompressible elastic fluid membrane in
three-dimensional space. We refer to \cite{Hu-Zhang} for more details.

\subsection{Geometric tensors and their evolution equations}
For a surface membrane $\Gamma=\RR(\vec{u},t)$ with a curve coordinate
$\vec{u}=(u^1,u^2)$, we can get the Frenet coordinate system
of the surface. Namely, the tangent vectors $\aaa_{\alpha}$ and the
unit normal vector $\nn$ are given by
$$\aaa_{\alpha}=\frac{\partial\RR}{\partial u^{\alpha}}\,(\al=1,2),\quad\nn=\frac{\aaa_1\times\aaa_2}{|\aaa_1\times\aaa_2|}.$$
The covariant metric tensor $(a_{\alpha\beta})_{1\le \al,\be\le 2}$
is defined as
$$a_{\alpha\beta}=\aaa_{\alpha}\cdot\aaa_{\beta}.$$
We denote its inverse by $(a^{\alpha\beta})_{1\le \al,\be\le 2}$,
which can be used to raise or lower the indices of the vectors and
tensors. For example,
$$
b_{\beta}^{\gamma}=a^{\alpha\gamma}b_{\alpha\beta}.
$$
The surface Christoffel symbols
$\Gamma_{\alpha\beta}^{\gamma}=\Gamma_{\beta\alpha}^{\gamma}$ and the
curvature tensor $b_{\alpha\beta}=b_{\beta\alpha}$ are given by the
Gauss-Weingarten-Codazzi equation:
\beno
&&\frac{\partial\aaa_{\alpha}}{\partial{u^{\beta}}}=
\Gamma_{\alpha\beta}^{\gamma}\aaa_{\gamma}+b_{\alpha\beta}\nn,\\
&&\frac{\partial\nn}{\partial{u^{\beta}}}=-b_{\beta}^{\gamma}\aaa_{\gamma}
=-a^{\alpha\gamma}b_{\alpha\beta}\aaa_{\gamma},\\
&&b_{\alpha\beta,\gamma}=b_{\alpha\gamma,\beta}.
\eeno
Here we use a comma followed by a lowercase Greek subscript to denote the covariant derivatives
based on the metric tensor $a_{\alpha\beta}$, that is,
\begin{equation}\label{eq:covariant derivative}
Q^{..\alpha.}_{.\beta..,\gamma}=\frac{\partial Q^{..\alpha.}_{.\beta..}}{\partial u^{\gamma}}
+\sum\Gamma_{\mu\gamma}^{\alpha}Q^{..\mu.}_{.\beta..}
-\sum\Gamma_{\beta\gamma}^{\mu}Q^{..\alpha.}_{.\mu..}.
\end{equation}
For example, we have
$$b_{\alpha\beta,\gamma}=\frac{\partial b_{\alpha\beta}}{\partial u^{\gamma}}
-\Gamma_{\alpha\gamma}^{\delta}b_{\delta\beta}
-\Gamma_{\beta\gamma}^{\delta}b_{\alpha\delta}.$$ Thus we can
rewrite the Gauss-Weigarten-Codarzzi equation as
\ben\label{eq:Guass-Weiggarten-Codazzi}
\aaa_{\alpha,\beta}=b_{\alpha\beta}\nn,~~
\nn_{,\alpha}=-b_{\alpha}^{\beta}\aaa_{\beta},~~
b_{\alpha\beta,\gamma}=b_{\alpha\gamma,\beta}. \een The mean
curvature $H$ and the Gaussian curvature $K$ of the surface are
given by \beno H=\frac{1}{2}b^{\alpha}_{\alpha},\quad
K=\frac{1}{2}(4H^2-b^{\alpha}_{\beta}b^{\beta}_{\alpha}). \eeno

In the following, let us derive the evolution equations of the geometric tensors. For this purpose,
we denote by $\vv(\vcu,t)$ the velocity of the surface given by
\ben\label{eq:surface velocity}
\vv(\vcu,t)=\frac{\partial{\RR(\vcu,t)}}{\partial t},
\een
and we decompose it into
\beno
\vv=v^{\alpha}\aaa_{\alpha}+v^{n}\nn.
\eeno
Using (\ref{eq:Guass-Weiggarten-Codazzi}), it is easy to find that
\begin{eqnarray}\label{eq:tangent-evolution}
\frac{\partial \aaa_{\alpha}}{\partial{t}}
&=&\frac{\partial}{\partial{t}}\frac{\partial{\RR}}{\partial u^{\alpha}}
=\frac{\partial}{\partial{u^{\alpha}}}\frac{\partial{\RR}}{\partial t}
=\frac{\partial\vv}{\partial{u^{\alpha}}}=\vv_{,\alpha}\nonumber\\
&=&(v^{\gamma}_{,\alpha}-v^{n}b_{\alpha}^{\gamma})\aaa_{\gamma}
+(v^{n}_{,\alpha}+v^{\gamma}b_{\alpha\gamma})\nn.
\end{eqnarray}
As $\nn\cdot\aaa_\al=0$, we get by (\ref{eq:tangent-evolution}) that
\beno
\frac{\partial{\nn}}{\partial{t}}\cdot\aaa_{\alpha}=-\frac{\partial
{\aaa_{\alpha}}}{\partial t}\cdot\nn=-(v^{n}_{,\alpha}+v^{\gamma}b_{\alpha\gamma}),
\eeno
which together with the fact of $\frac{\partial{\nn}}{\partial{t}}\cdot\nn=0$ implies that
\ben\label{eq:normal evolution}
\frac{\partial{\nn}}{\partial{t}}=-(v^{\beta}b^{\alpha}_{\beta}+a^{\al\be}v^n_{,\be})\aaa_\al.
\een
The evolution equation of the metric tensor is given by
\begin{eqnarray}\label{eq:metric evolution}
\frac{\partial{a_{\al\be}}}{\partial{t}}=
\frac{\partial{\aaa_{\alpha}}}{\partial{t}}\cdot\aaa_{\beta}+\frac{\partial
{\aaa_{\beta}}}{\partial t}\cdot\aaa_{\alpha}
=(v_{\alpha,\beta}+v_{\beta,\alpha})-2v^nb_{\al\be}.
\end{eqnarray}
Differentiating the identity $a^{\alpha\beta}a_{\beta\gamma}=\delta_{\gamma}^{\alpha}$ with respect to $t$,
we get by (\ref{eq:metric evolution}) that
\begin{eqnarray}\label{eq:covariant metric evolution}
\frac{\partial a^{\alpha\beta}}{\partial{t}}=-a^{\alpha\gamma}a^{\beta\delta}(v_{\gamma,\delta}+v_{\delta,\gamma})+2v^{n}b^{\alpha\beta}.
\end{eqnarray}
And differentiating $b_{\alpha\beta}=-\aaa_{\alpha}\cdot\nn_{,\beta},$ we get by (\ref{eq:Guass-Weiggarten-Codazzi})--(\ref{eq:normal evolution}) that
\begin{eqnarray}\label{eq:curvature evolution}
\frac{\partial{b_{\al\be}}}{\partial{t}}=(v^n_{,\al\be}-v^nb^{\gamma}_{\alpha}b_{\gamma\beta})+(v^{\gamma}_{,\beta}b_{\alpha\gamma}
+v^{\gamma}_{,\alpha}b_{\gamma\beta})+v^{\gamma}b_{\al\be,\gamma}.
\end{eqnarray}
Due to $2H=a^{\alpha\beta}b_{\alpha\beta}$, we get by (\ref{eq:covariant metric evolution}) and (\ref{eq:curvature evolution}) that
\begin{eqnarray}\label{eq:mean curvature evolution}
2\frac{\partial{H}}{\partial{t}}=a^{\al\be}v^n_{,\al\be}+v^nb^{\alpha}_{\beta}b^{\beta}_{\alpha}+2v^{\alpha}H_{,\alpha}.
\end{eqnarray}

\subsection{The derivation of the elastic surface model}

In this subsection, we choose $\vec{u}=(u^1, u^2)$ as the Lagrangian coordinate of the
moving fluid surface. In this coordinate system, the velocity of the fluid on the surface is equal
to the velocity of the fluid $\vv$ given by (\ref{eq:surface velocity}).
The Helfrich bending elastic energy {\cite{DHelfrich, Helfrich}} is
\begin{equation}\label{eq:energy-4}
E_H=\int_{\Gamma} C_1^{\alpha\beta\gamma\delta}(B_{\alpha\beta}-
b_{\alpha\beta})(B_{\gamma\delta}-b_{\gamma\delta}) \ud S,
\end{equation}
where $B_{\alpha\beta}$ is the spontaneous curvature tensor,
and the fourth-order tensor $C_1^{\alpha\beta\gamma\delta}$ is given by
\beno
C_1^{\alpha\beta\gamma\delta}=(k_1-\varepsilon_1)a^{\alpha\beta}a^{\gamma
\delta}+\varepsilon_1(a^{\alpha\gamma}a^{\beta\delta}+a^{\alpha\delta}
a^{\beta\gamma}),
\eeno
where $k_1$ and $\varepsilon_1$ are positive elastic coefficients and $k_1\geq\varepsilon_1$.

As the membrane is a two-dimensional incompressible fluid, we have
$$\nabla_{\Gamma}\cdot\vv=0,$$
which is equivalent to
\begin{equation}\label{eq:incompressible}
v^{\alpha}_{,\alpha}=2Hv^{n}.
\end{equation}
Then by applying the principle of virtual work, we obtain
elastic stresses. For isotropic Newtonian membrane fluids, the dynamical
equation  of the membrane is
\begin{eqnarray}\label{eq:velocity}
&&\varrho\frac{\partial\vv}{\partial t}=\mathbf{f}+(T^{\alpha\beta}\aaa_{\beta})_{,\alpha}
+(q^{\alpha}\nn)_{,\alpha},
\end{eqnarray}
where $\varrho$ is the membrane fluid density. The in-plane stress
tensor $T^{\alpha\beta}$ and transverse shear stress $q^{\alpha}$
are given by
\beno T^{\alpha\beta}&=&-\Pi
a^{\alpha\beta}+J^{\alpha\beta}+M^{\alpha\mu}b_{\mu}^{\beta},
\\
q^{\alpha}&=&M^{\alpha\beta}_{,\beta},\\
J^{\alpha\beta}&=&C^{\alpha\beta\gamma\delta}S_{\gamma\delta},\\
C^{\alpha\beta\gamma\delta}&=&(k_0-\varepsilon_0)a^{\alpha\beta}a^{\gamma
\delta}+\varepsilon_0(a^{\alpha\gamma}a^{\beta\delta}+a^{\alpha\delta}a^{\beta\gamma}),\\
M^{\alpha\beta}&=&C_1^{\alpha\beta\gamma\delta}(B_{\gamma\delta}-b_{\gamma\delta}),
\eeno
where $\Pi$ is the surface pressure (tension), and the rate of the surface strain is given by
$$S_{\alpha\beta}\eqdef\frac{1}{2}\frac{\partial{a_{\ab}}}{\partial{t}}=\frac{1}{2}(v_{\alpha,\beta}
+v_{\beta,\alpha})-v^{n}b_{\alpha\beta}.$$
From (\ref{eq:incompressible}), it is easy to see that $S^\al_\al=0$, and thus $J^{\al\be}=2\varepsilon_0S^{\al\be}$.
Furthermore, the above equations have the following energy
dissipation relation:
\ben\label{eq:energy-s2}
\frac{1}{2}\frac{\ud}{\ud t}\Big(E_H+\int_{\Gamma}\varrho|\vv|^2\ud
S\Big)=
-\int_{\Gamma}C^{\alpha\beta\gamma\delta}S_{\alpha\beta}S_{\gamma\delta}\ud
S =-2\varepsilon_0\int_{\Gamma}S^{\al\be}S_{\alpha\beta}\ud S. \een
If the function $B(\vcu)=B(\vcu)a_{\alpha\beta}$ with $B$ independent of the time $t$, then (\ref{eq:energy-4})
can be reduced to
\beno E_H=\int_{\Gamma}4k_1(H-B)^2+4\varepsilon_1(H^2-K)\ud{S},
\eeno
and the velocity equation (\ref{eq:velocity}) can be reduced to the following form
\ben\label{eq:velocity-2}
\varrho\frac{\partial{\vv}}{\partial{t}}&=&(-\Pi{a^{\alpha\beta}\aaa_{\alpha}})_{,\beta}+2\varepsilon_0
(S^{\al\be}\aaa_{\alpha})_{,\beta}-4k_1Ha^{\al\be}B_{,\beta}\aaa_{\alpha}\nonumber\\
&&+2k_1(\Delta_{\Gamma}{B}-2BK)\nn-2(k_1+\varepsilon_1)\Big(\Delta_{\Gamma}{H}
+2H(H^2-K)\Big)\nn. \een Here $\Delta_\Gamma$ denotes the
Lapalace-Beltrami operator on the surface $\Gamma$, and  $K$ is the
Gaussian curvature. We refer to the appendix for the derivations of
(\ref{eq:energy-s2}) and (\ref{eq:velocity-2}).

In this paper, we only consider the simple case with $B(\vcu)\equiv 0$.
By the rescaling argument, we can set $\varrho=1, 4(k_1+\varepsilon_1)=1$. Thus, we obtain
\begin{eqnarray}\label{eq:model-lagrangian}
\left\{
\begin{array}{l}
\frac{\partial{\vv}}{\partial{t}}=(-\Pi{a^{\alpha\beta}\aaa_{\alpha}})_{,\beta}+2\varepsilon_0
(S^{\al\be}\aaa_{\alpha})_{,\beta}
-\frac{1}{2}\Big(\Delta_{\Gamma}{H}+2H(H^2-K)\Big)\nn,\\
v^{\alpha}_{,\alpha}=2Hv^{n}.
\end{array}\right.
\end{eqnarray}

\setcounter{equation}{0}
\section{New formulation of the system}

Motivated by \cite{Ambrose}, we will reformulate (\ref{eq:model-lagrangian}) into a new system in the isothermal coordinates.
That is, we choose a coordinate $\vec{x}=(x^1, x^2)$ such that the metric tensor satisfies
\begin{equation}\label{eq:isotherma relation}
a_{11}=a_{22},\quad a_{12}=a_{21}=0.
\end{equation}
As the tangential velocity of the surface only serves to reparameterize the surface,
in the following we choose them such that
\begin{equation}\label{eq:isothermal relation-2}
\frac{\partial a_{11}}{\partial t}=\frac{\partial a_{22}}{\partial t},\quad
\frac{\partial a_{12}}{\partial t}= 0.
\end{equation}
As a result, if it holds for the initial surface
the relation (\ref{eq:isotherma relation}) will be preserved for any time $t$.

\subsection{Elliptic system for the tangential velocity of the surface}

In the sequel, for convenience, we denote $f_{\alpha}=\frac{\partial f}{\partial
x^1}, f_{\beta}=\frac{\partial f}{\partial x^2}$ , and whereas
 $(\cdot)_{,\alpha}$ (or $(\cdot)_{,\beta}$) denotes the covariant derivative with respect to $x^1$ (or $x^2$).
 The unit tangent vector and the unit normal
vector of the surface are given respectively by
\begin{equation}
\ta=\frac{\RR_{\alpha}}{|\RR_{\alpha}|},\quad \tb=\frac{\RR_{\beta}}{|\RR_{\beta}|},\quad
\nn=\frac{\RR_{\alpha}\times\RR_{\beta}}{|\RR_{\alpha}\times\RR_{\beta}|}.
\end{equation}
We denote $E, F$, and $G$ by the coefficients of the first fundamental form, and $L, M,$ and $N$ by
the coefficients of the second fundamental form. Namely,
\beno
&&E=\RR_\al\cdot\RR_\al,\quad F=\RR_\al\cdot\RR_\be,\quad G=\RR_\be\cdot\RR_\be,\\
&&L=\RR_{\al\al}\cdot\nn,\quad M=\RR_{\al\be}\cdot\nn,\quad N=\RR_{\be\be}\cdot\nn.
\eeno
In the isothermal coordinates, we have $E=G, F=0$. The Christoffel symbols can be calculated as follows:
\ben\label{Christoffel symbol}
\begin{array}{l}
\Gamma_{11}^1=\Gamma_{12}^2=\Gamma_{21}^2=-\Gamma_{22}^1=\frac{\Ea}{2E},\\
\Gamma_{12}^1=\Gamma_{21}^1=\Gamma_{22}^2=-\Gamma_{11}^2=\frac{\Eb}{2E}.
\end{array}
\een
And the following identities can be verified easily:
\begin{eqnarray}\label{eq:isorelations}
\begin{array}{l}
\taa\cdot\tb=-\tba\cdot\ta=-\Big(\frac{\RRb}{|\RRb|}\Big)_{\alpha}\cdot\frac{\RRa}{|\RRa|}=-\frac{\RR_{\ab}\cdot\RRa}{E}=-\frac{\Eb}{2E},\\
\tbb\cdot\ta=-\tab\cdot\tb=-\Big(\frac{\RRa}{|\RRa|}\Big)_{\beta}\cdot\frac{\RRb}{|\RRb|}=-\frac{\RR_{\ab}\cdot\RRb}{E}=-\frac{\Ea}{2E},\\
\taa\cdot\nn=-\ta\cdot\nna=-\frac{\RRa}{|\RRa|}\cdot\nna=\frac{L}{\se},\\
\tba\cdot\nn=-\tb\cdot\nna=-\frac{\RRb}{|\RRb|}\cdot\nna=\frac{M}{\se},\\
\tab\cdot\nn=-\ta\cdot\nnb=-\frac{\RRa}{|\RRa|}\cdot\nnb=\frac{M}{\se},\\
\tbb\cdot\nn=-\tb\cdot\nnb=-\frac{\RRb}{|\RRb|}\cdot\nnb=\frac{N}{\se}.
\end{array}
\end{eqnarray}

For a given normal velocity $\un(x^1,x^2,t)$, we assume that the evolution of the surface is determined by
\begin{equation}\label{eq:velocity of surface}
\frac{\partial\RR(\vec{x}, t)}{\partial
t}=\un\nn+W_1\ta+W_2\tb\eqdef\ww.
\end{equation}
Then it follows from (\ref{eq:isorelations}) that
\begin{eqnarray}
\RR_{\alpha t}&=&(\una+\frac{W_1L}{\se}+\frac{W_2M}{\se})\nn+(W_{1\alpha}-\frac{\un L}{\se}
+\frac{W_2\Eb}{2E})\ta\nonumber\\&&
+(W_{2\alpha}-\frac{\un M}{\se}-\frac{W_1\Eb}{2E})\tb\nonumber\\
&\eqdef&A_1\nn+A_{01}\ta+A_3\tb,\label{eq:R-evolution-1}\\
\RR_{\beta t}&=&(\unb+\frac{W_1M}{\se}+\frac{W_2N}{\se})\nn+(W_{1\beta}-\frac{\un M}{\se}
-\frac{W_2\Ea}{2E})\ta
\nonumber\\&&
+(W_{2\beta}-\frac{\un N}{\se}+\frac{W_1\Ea}{2E})\tb\nonumber\\
&\eqdef&A_2\nn+A_4\ta+A_{02}\tb.\label{eq:R-evolution-2}
\end{eqnarray}
Consequently, we obtain
\beno
&&E_t=2\RR_{\alpha t}\cdot\RRa=2\se\RR_{\alpha t}\cdot\ta=2\se A_{01}, \\
&&G_t=2\RR_{\beta  t}\cdot\RRb=2\se\RR_{\beta  t}\cdot\tb=2\se A_{02},\\
&&F_t=\se(\RR_{\alpha t}\cdot\tb+\RR_{\beta t}\cdot\ta)=\se(A_3+A_4).
\eeno
Now the relation (\ref{eq:isothermal relation-2}) is equivalent to
\beno
(E-G)_t=0,\quad F_t=0,
\eeno
which implies that
\begin{eqnarray}\label{eq:tangential-elliptic system}
\left\{
\begin{array}{l}
\Big(\frac{W_1}{\se}\Big)_{\alpha}-\Big(\frac{W_2}{\se}\Big)_{\beta}=\frac{\un(L-N)}{E},\\
\Big(\frac{W_1}{\se}\Big)_{\beta}+\Big(\frac{W_2}{\se}\Big)_{\alpha}=\frac{2\un M}{E}.
\end{array}\right.
\end{eqnarray}
This is an elliptic system for $(W_1, W_2)$. As mentioned above, if the surface evolves as (\ref{eq:velocity of surface})
with $(W_1, W_2)$ determined by (\ref{eq:tangential-elliptic system}), the coordinate will always be isothermal.

\begin{Remark}
The above system can also be obtained by using (\ref{eq:metric evolution}) directly.
\end{Remark}

Let us conclude this section by deriving the elliptic equations for $E$ and $\RR$.
Noticing that $E=\RRa\cdot\RRa=\RRb\cdot\RRb$ and $\RRa\cdot\RRb=0$, we have
\begin{eqnarray}\label{eq:E-elliptic equation}
\Delta E&=&(\RRa\cdot\RRa)_{\beta\beta}+(\RRb\cdot\RRb)_{\alpha\alpha}-2(\RRa\cdot\RRb)_{\alpha\beta}\nonumber\\
&=&2(\RR_{\alpha\beta}\cdot\RR_{\alpha\beta}-\RR_{\alpha\alpha}\cdot\RR_{\beta\beta}).
\end{eqnarray}
On the other hand, we have
\beno
\RR_{\alpha}\cdot\Delta\RR=0,\quad \nn\cdot\Delta\RR=L+N=2EH,
\eeno
which means that
\ben\label{eq:R-elliptic equation}
\Delta\RR=2EH\nn.
\een

\begin{Remark}
From (\ref{eq:E-elliptic equation}) and the standard elliptic estimate,
it is easy to find that $E$ has the same
regularity as the surface. This fact is  noted by S.-S
Chern in \cite{Chern}. Then we can gain two more regularities of  $\RR$ from the regularity of the mean curvature $H$ by
using (\ref{eq:R-elliptic equation}).
Specifically, we will use (\ref{eq:E-elliptic equation}) and (\ref{eq:R-elliptic equation})
to construct our approximate solutions in Section 5.1.

\end{Remark}

\subsection{The velocity equation in the isothermal coordinate}

Assume that $\uu(\vcx,t)$ is the velocity of the fluid in the isothermal coordinate. Hence,
$\vv(\vcu,t)=\uu(\vcx(\vcu,t),t)$, where $\vcu$ is the Lagrangian coordinate.
And, we have
\beno
\uu(\vcx(\vcu,t),t)&=&\vv(\vcu,t)=\frac{\ud\RR(\vcx(\vcu,t),t)}{\ud t}\nonumber\\
&=&\frac{\partial{\RR}}{\partial t}\circ\vcx
+\frac{\partial{x^1({\vcu,t})}}{\partial{t}}\frac{\partial{\RR}}{\partial x^1}\circ\vcx
+\frac{\partial{x^2({\vcu,t})}}{\partial{t}}\frac{\partial{\RR}}{\partial x^2}\circ\vcx\nonumber\\
&=&\ww\circ\vcx+x^1_t\RRa\circ\vcx+x^2_t\RRb\circ\vcx.
\eeno
Hence, we have
$$x^1_t(\vcu,t)=\big[\frac{1}{\se}(\uu-\ww)\cdot\ta\big]\circ\vcx,\quad
x^2_t(\vcu,t)=\big[\frac{1}{\se}(\uu-\ww)\cdot\tb\big]\circ\vcx.$$
Consequently,
\begin{eqnarray}
\frac{\partial\vv(\vcu,t)}{\partial{t}}&=&\frac{\ud\uu(\vcx(\vcu,t),t)}{\ud t}\nonumber\\
&=&\frac{\partial{\uu}}{\partial t}\circ\vcx
+\frac{\partial{x^1({\vcu,t})}}{\partial{t}}\frac{\partial{\uu}}{\partial x^1}\circ\vcx
+\frac{\partial{x^2({\vcu,t})}}{\partial{t}}\frac{\partial{\uu}}{\partial x^2}\circ\vcx\nonumber\\
&=&\Big(\frac{\partial{\uu}}{\partial t}+\frac{1}{\se}\big[(\uu-\ww)\cdot\ta\big]\uupa
+\frac{1}{\se}\big[(\uu-\ww)\cdot\tb\big]\uupb\Big)\circ\vcx.\nonumber
\end{eqnarray}
The above equation can also be derived by Oldroyd's theorem \cite{Oldroyd}.

On the other hand,  by (\ref{eq:Guass-Weiggarten-Codazzi}) we have that
\beno
&&\RR_{\al,\al}=L\nn,\quad  \RR_{\al,\be}=M\nn,\quad \RR_{\be,\be}=N\nn,\\
&&H=\f 1{2E}(L+N),\quad \Delta_\Gamma H=\f 1 E\Delta H.
\eeno
Given that the right-hand side of (\ref{eq:model-lagrangian}) is coordinate-invariant,
the first equation of (\ref{eq:model-lagrangian}) can be reduced to
\begin{eqnarray}\label{eq:velocity equation-isothermal}
&&\frac{\partial{\uu}}{\partial{t}}+\frac{1}{\se}[(\uu-\ww)\cdot\ta]\uu_{\alpha}+\frac{1}{\se}[(\uu-\ww)\cdot\tb]\uu_{\beta}\nonumber\\
&&=-2H\Pi\nn-\frac{\Pi_{\alpha}}{\se}\ta-\frac{\Pi_{\beta}}{\se}\tb+\frac{2\epsilon_0}{E^2}\Big(S_{11}L+2S_{12}M+S_{22}N\Big)\nn\nonumber\\
&&\quad+{2\epsilon_0}\se(S^{11}_{,\alpha}+S^{12}_{,\beta})\ta
+{2\epsilon_0}{\se}(S^{12}_{,\alpha}+S^{22}_{,\beta})\tb
\nonumber\\
&&\quad-\frac{1}{2}\Big(\frac{\Delta H}{E}+\frac{H}{2E^2}\big((L-N)^2+4M^2\big)\Big)\nn.
\end{eqnarray}
Here $S_{11}=\uu_\al\cdot \RR_\al, S_{12}=\f12(\uu_\al\cdot \RR_\be+\uu_\be\cdot\RR_\al),
S_{22}=\uu_\be\cdot \RR_\be$. So,
\beno
&&S_1^1=\frac{S_{11}}{E},\quad S_1^2=S_2^1=\frac{S_{12}}{E},\quad S_2^2=\frac{S_{22}}{E},\\
&&S^{11}=\frac{S_{11}}{E^2},\quad S^{12}=S^{21}=\frac{S_{12}}{E^2},\quad S^{22}=\frac{S_{22}}{E^2}.
\eeno

\subsection{New equivalent system}

Setting $\uu(\vec{x},t)=\un\nn+\uua\ta+\uub\tb$, we infer from (\ref{eq:isorelations}) that
\begin{eqnarray}
\uupa&=&(\una+\frac{L\uua+M\uub}{\se})\nn+(\uaa-\frac{L\un}{\se}+\frac{\uub\Eb}{2E})\ta\nonumber\\
&&+(\uba-\frac{M\un}{\se}-\frac{\uua\Eb}{2E})\tb\nonumber\\
&\eqdef& f_1\nn+g_{11}\ta+g_{12}\tb,\label{eq:u_al}\\
\uupb&=&(\unb+\frac{M\uua+N\uub}{\se})\nn+(\uab-\frac{M\un}{\se}-\frac{\uub\Ea}{2E})\ta\nonumber\\&&
+(\ubb-\frac{N\un}{\se}+\frac{\uua\Ea}{2E})\tb\nonumber \\
&\eqdef& f_2\nn+g_{21}\ta+g_{22}\tb.\label{eq:u_be}
\end{eqnarray}
Using (\ref{eq:R-evolution-1}) and (\ref{eq:R-evolution-2}), we find that
\beno
&&\f {\p\ta}{\p t}=\f 1 {\sqrt{E}}(A_1\nn+A_3\tb),\\
&&\f {\p\tb}{\p t}=\f 1 {\sqrt{E}}(A_2\nn+A_4\ta),\\
&&\f {\p\nn}{\p t}=-\f 1 {\sqrt{E}}(A_1\ta+A_2\tb).
\eeno
Thus, the equation (\ref{eq:velocity equation-isothermal}) can be rewritten as
\begin{eqnarray}\label{eq:velocity equation-isothermal-component}
&&\Big(\un_t+\frac{1}{\se}(\uua A_1+\uub{
A_2})+\frac{\uua-W_1}{\se}f_1+\frac{\uub-W_2}{\se}f_2\Big)\nn\nonumber\\
&&\quad+\Big(U_{1t}+\frac{1}{\se}(\uub A_4-\un
{A_1})+\frac{\uua-W_1}{\se}g_{11}+\frac{\uub-W_2}{\se}g_{21}\Big)\ta
\nonumber\\
&&\quad+\Big(U_{2t}+\frac{1}{\se}(\uua A_3-\un
{A_2})+\frac{\uua-W_1}{\se}g_{12}+\frac{\uub-W_2}{\se}g_{22}\Big)\tb
\nonumber\\
&&=-2H\Pi\nn-\frac{\Pi_{\alpha}}{\se}\ta-\frac{\Pi_{\beta}}{\se}\tb+\frac{2\epsilon_0}{E^2}\Big(S_{11}L+2S_{12}M+S_{22}N\Big)\nn\nonumber\\
&&\quad+\frac{2\epsilon_0}{E\se}(S_{11,\alpha}+S_{12,\beta})\ta
+\frac{2\epsilon_0}{E\se}(S_{12,\alpha}+S_{22,\beta})\tb \nonumber\\&&
\quad-\frac{1}{2}\Big(\frac{\Delta
H}{E}+\frac{H}{2E^2}\big((L-N)^2+4M^2\big)\Big)\nn.
\end{eqnarray}
Here we used the fact that $E^2\big(S^{11}_{,\alpha}+S^{21}_{,\beta})=S_{11,\alpha}+S_{12,\beta}$.

Now we calculate $S_{11,\alpha}+S_{12,\beta}$.
By (\ref{eq:u_al}) and (\ref{eq:u_be}), the surface strain rate tensor can be written as
\begin{eqnarray}
S_{11}&=&\uupa\cdot\RRa=\se\uaa+\frac{\uub\Eb}{2\se}-L\un\nonumber\\
&=&(\se\uua)_{\alpha}+\frac{\Eb\uub-\Ea\uua}{2\se}-L\un,\label{eq:S_11}\\
S_{12}&=&S_{21}=\frac{1}{2}\big(\uupa\cdot\RRb+\uupb\cdot\RRa\big)\nonumber\\
&=&\frac{\se}{2}\big(\uab+\uba\big)-\frac{\Ea\uub+\Eb\uua}{4\se}-M\un\nonumber\\
&=&\frac{1}{2}(\se\uua)_{\beta}+\frac{1}{2}(\se\uub)_{\alpha}-\frac{\Ea\uub+\Eb\uua}{2\se}-M\un,\label{eq:S_12}\\
S_{22}&=&\uupb\cdot\RRb=\se\ubb+\frac{\Ea\uua}{2\se}-N\un\nonumber\\
&=&(\se\uub)_{\beta}+\frac{\Ea\uua-\Eb\uub}{2\se}-N\un.\label{eq:S_22}
\end{eqnarray}
Since the incompressible condition $\nabla_{\Gamma}\cdot\uu=0$ can also be written as
$S_{11}+S_{22}=0$, we get by (\ref{eq:S_11}) and (\ref{eq:S_22}) that
\begin{equation}\label{eq:velocity free}
(\se\uua)_{\alpha}+(\se\uub)_{\beta}-2EH\un=0.
\end{equation}
We get by (\ref{Christoffel symbol}) and $S_{11}+S_{22}=0$ that
\begin{eqnarray}
S_{11,\alpha}+S_{12,\beta}&=&
S_{11\alpha}-2\Gamma_{11}^{i}S_{1i}+S_{12\beta}-\Gamma_{22}^iS_{1i}-\Gamma_{12}^iS_{i2} \nonumber\\
&=& S_{11\alpha}+S_{12\beta}.\nonumber
\end{eqnarray}
Similarly,
\beno
S_{12,\alpha}+S_{22,\beta}=S_{12\alpha}+S_{22\beta}.
\eeno
Using (\ref{eq:S_11})-(\ref{eq:velocity free}), we find that
\begin{eqnarray}
S_{11\alpha}+S_{12\beta}&=&\Big((\se\uua)_{\alpha}+\frac{\Eb\uub-\Ea\uua}{2\se}-L\un\Big)_{\alpha}\nonumber\\
&&+\bigg(\frac{1}{2}(\se\uua)_{\beta}+\frac{1}{2}(\se\uub)_{\alpha}-\frac{\Ea\uub+\Eb\uua}{2\se}-M\un\bigg)_{\beta}\nonumber\\
&=&\frac{\Delta(\se\uua)}{2}+\frac{1}{2}\Big((\se\uua)_{\alpha}+(\se\uub)_{\beta}\Big)_{\alpha}  \nonumber\\
&&+\Big(\frac{\Eb\uub-\Ea\uua}{2\se}-L\un\Big)_{\alpha}-\Big(\frac{\Ea\uub+\Eb\uua}{2\se}+M\un\Big)_{\beta}\nonumber\\
&=&\frac{\Delta(\se\uua)}{2}+\Big(\frac{\Eb\uub-\Ea\uua}{2\se}-L\un+EH\un\Big)_{\alpha}\nonumber\\&&
-\Big(\frac{\Ea\uub+\Eb\uua}{2\se}+M\un\Big)_{\beta}.\nonumber
\end{eqnarray}
Thus, we obtain  from (\ref{eq:velocity equation-isothermal-component}) the evolution equation for $U_1$:
\ben\label{eq:U_1}
\frac{\partial \uua}{\partial
t}&=&\frac{\epsilon_0}{E\se}\Delta(\se\uua)-\frac{\Pi_{\alpha}}{\se}\nonumber\\
&&+\frac{2\epsilon_0}{E\se}\bigg(\frac{\Eb\uub-\Ea\uua}{2\se}-L\un+EH\un\bigg)_{\alpha}\nonumber\\
&&+\frac{1}{2\se}\Big((\un)^2\Big)_{\alpha}-\frac{\epsilon_0}{E\se}\bigg(\frac{\Ea\uub+\Eb\uua}{\se}+2M\un\bigg)_{\beta}\nonumber\\
&&-\frac{1}{\se}\Big(\uub W_{1\beta}+(\uua-\Wa)\uaa+(\uub-\Wb)\uab\Big)\nonumber\\
&&+\frac{1}{E}(2\uub\un M+\uua\un L)-\frac{1}{2E\se}\Big((\uua-\Wa)\uub\Eb-\uub^2\Ea\Big).
\een
A similar evolution equation for $U_2$ can also be obtained, although we omit the details here.
The evolution equation for the normal velocity $U^n$ is
\ben\label{eq:U_n}
\frac{\partial \un}{\partial t}&=&-\frac{1}{2E}\Delta
H-2H\Pi-\frac{1}{\se}\Big((2\uua-\Wa)\un_{\alpha}+(2\uub-W_2)\un_{\beta}\Big)\nonumber\\
&&+\frac{2\epsilon_0}{E\se}(L\uaa+M\uab+M\uba+N\ubb)\nonumber\\
&&+\frac{\epsilon_0}{E^2\se}(L\Eb\uub-M\Ea\uub-M\Eb\uua+N\Ea\uua)\nonumber\\
&&-\frac{2\epsilon_0U^n}{E\se}(L^2+2M^2+N^2)\nonumber\\
&&-\frac{H}{4E^2}\Big((L-N)^2+4M^2\Big)-\frac{1}{E}(L\uua^2+2M\uua\uub+N\uub^2).
\een
Due to (\ref{eq:mean curvature evolution}), the evolution equation for the mean curvature $H$ is
\begin{eqnarray}\label{eq:H}
\frac{\partial{H}}{\partial{t}}=\frac{1}{2E}\Delta\un
+\frac{W_1}{\se}H_{\alpha}+\frac{W_2}{\se} H_{\beta}+\frac{\un}{2E^2}(L^2+2M^2+N^2).
\end{eqnarray}
We denote
\beno
&&F_1=\frac{\Eb\uub-\Ea\uua}{2\se}+\frac{\un}{2}(N-L),\quad F_2=\frac{1}{2}(\un)^2,\\
&&F_3=\frac{\Ea\uub+\Eb\uua}{\se}+2M\un,\\
&&F_4=-\frac{1}{\se}\uub W_{1\beta}+\frac{1}{E}(2M\uub\un+L\uua\un)-\frac{1}{2E\se}\Big((\uua-\Wa)\uub\Eb-\uub^2\Ea\Big),\\
&&F_5=-\frac{1}{\se}\uua W_{2\alpha}+\frac{1}{E}(2M\uua\un+N\uub\un)-\frac{1}{2E\se}\Big((\uub-\Wb)\uua\Ea-\uua^2\Eb\Big),\\
&&F_6=\frac{\epsilon_0}{E^2\se}(L\Eb\uub-M\Ea\uub-M\Eb\uua+N\Ea\uua)
-\frac{2\epsilon_0U^n}{E\se}(L^2+2M^2+N^2)\\
&&\qquad-\frac{H}{4E^2}\Big((L-N)^2+4M^2\Big)-\frac{1}{E}(L\uua^2+2M\uua\uub+N\uub^2),\\
&&F_7=\frac{\un}{2E^2}(L^2+2M^2+N^2),
\eeno
and $u_i=\frac{U_i}{\se},\, w_i=\frac{W_i}{\se}$ for $i=1,2$.

From (\ref{eq:velocity free})-(\ref{eq:H}), we obtain the following equivalent full system:
\begin{eqnarray}
&&\frac{\partial\RR}{\partial t}=\un\nn+W_1\ta+W_2\tb,\label{eq:full system-R}\\
&&\frac{\partial \uua}{\partial t}=\epsilon_0\f 1 {E\sqrt{E}}\Delta(\sqrt{E}\uua)-\f 1 {\sqrt{E}}{\Pi_{\alpha}}
+\f {2\epsilon_0}{E\sqrt{E}}F_{1\alpha}+\frac 1{2\sqrt{E}}F_{2\alpha}-\f {\epsilon_0}{E\sqrt{E}}F_{3\beta}\nonumber\\
&&\qquad\quad-(u_1-w_1)\uaa-(u_2-w_2)\uab+F_4,\label{eq:full system-U1}\\
&&\frac{\partial \uub}{\partial t}=\epsilon_0\f 1 {E\sqrt{E}}\Delta(\sqrt{E}\uub)-\f 1 {\sqrt{E}}{\Pi_{\beta}}
-\f {2\epsilon_0}{E\sqrt{E}}F_{1\beta}+\frac{1}{2\sqrt{E}}F_{2\beta}-\f {\epsilon_0}{E\sqrt{E}}F_{3\alpha}\nonumber\\
&&\qquad\quad-(u_1-w_1)\uba-(u_2-w_2)\ubb+F_5,\label{eq:full system-U2}\\
&&\frac{\partial \un}{\partial t}=-\frac{1}{2E}\Delta
H-2H\Pi-\Big((2u_1-w_1)\un_{\alpha}+(2u_2-w_2)\un_{\beta}\Big)\nonumber\\
&&\qquad\quad+\f {2\epsilon_0}{E\sqrt{E}}(L\uaa+M\uab+M\uba+N\ubb)+F_6,\label{eq:full system-Un}\\
&&\frac{\partial{H}}{\partial{t}}=\frac{1}{2E}\Delta\un
+w_1H_{\alpha}+w_2 H_{\beta}+F_7,\label{eq:full system-H}\\
&&(\se\uua)_{\alpha}+(\se\uub)_{\beta}-2EH\un=0,\label{eq:full system-velocity free}
\end{eqnarray}
where $(W_1,W_2)$ is determined by the elliptic system (\ref{eq:tangential-elliptic system}).
\begin{Remark}\label{rem:pressure}
Actually, (\ref{eq:full system-H}) is induced by (\ref{eq:full
system-R}). However, in order to perform a suitable energy estimate,
we add it to the full system.
\end{Remark}

\subsection{The equation of the pressure}

Using the incompressible condition $(\se\uua)_{\alpha}+(\se\uub)_{\beta}=2EH\un$, we find that
\beno
(\se U_{1t})_{\alpha}+(\se U_{2t})_{\beta}-2EH\un_t=\un(2EH)_t-\Big(
\frac{E_t}{2\se}\uua\Big)_{\alpha}-\Big(\frac{E_t}{2\se}\uub\Big)_{\beta}.
\eeno
We denote the left-hand side of (\ref{eq:velocity equation-isothermal-component}) by
\beno
(U_t^n+\f {K^n} {\sqrt{E}})\nn+(U_{1t}+\f {K_1} {\sqrt{E}})\ta+(U_{2t}+\f {K_2} {\sqrt{E}})\tb.
\eeno
Noting that $E_t=2\sqrt{E}A_{01}$, then we have by (\ref{eq:velocity equation-isothermal-component}) that
\beno
&&-\Pi_{\alpha\alpha}+2\epsilon_0\Big(\frac{S_{11\alpha}+S_{21\beta}}{E}\Big)_{\alpha}
-\Pi_{\beta\beta}+2\epsilon_0\Big(\frac{S_{12\alpha}+S_{22\beta}}{E}\Big)_{\beta}+4EH^2\Pi\\
&&-K_{1\alpha}-K_{2\beta}+2\se HK^n-\frac{4\epsilon_0
H}{E}\big(S_{11}L+2S_{12}M+S_{22}N\big)\\&&
+H\Big(\Delta H+\frac{H}{2E}\big((L-N)^2+4M^2\big)\Big)\\
&&=\un(2EH)_t-(A_{01}\uua)_{\alpha}-(A_{01}\uub)_{\beta}.
\eeno
By a direct computation, we obtain
\begin{eqnarray*}
&&\un(2EH)_t+(K_1-A_{01}\uua)_{\alpha}+(K_2-A_{01}\uub)_{\beta}-2\sqrt{E}HK^n\\
&&=\Big(\uua\uaa+\uub\uab-\frac{\un}{\se}(2M\uub-(L-N)\uua)+\frac{\uua\uub{\Eb}-\uub^2\Ea}{2E}\Big)_{\alpha}\\
&&\quad+\Big(\uua\uba+\uub\ubb-\frac{\un}{\se}(2M\uua-(N-L)\uub)+\frac{\uua\uub{\Ea}-\uua^2\Eb}{2E}\Big)_{\beta}\\
&&\quad-4H\se(\uua\una+\uub\unb)-2H(L\uua^2+2M\uub\uub+N\uub^2)-(\una)^2-(\unb)^2\\
&&\quad+\frac{(\un)^2}{E}(L^2+2M^2+N^2-4E^2H^2)\eqdef\cG_1.
\end{eqnarray*}
And using the incompressible condition again, we get
\begin{eqnarray*}
&&\Big(\frac{S_{11\alpha}+S_{21\beta}}{E}\Big)_{\alpha}+
\Big(\frac{S_{12\alpha}+S_{22\beta}}{E}\Big)_{\beta}\nonumber\\
&&=-\frac{\Ea}{2E^2}(S_{11\alpha}+S_{21\beta})-\frac{\Eb}{2E^2}(S_{12\alpha}+S_{22\beta})+\f 1 E\Delta(2EH\un)\\
&&\quad+\f 1 E\Big((\se)_{\beta}\uub-(\se)_{\alpha}\uua-L\un\Big)_{\alpha\alpha}
-\f 2 E\Big((\se)_{\alpha}\uub+(\se)_{\beta}\uua+M\un\Big)_{\alpha\beta}\\
&&\quad+\f 1 E\Big((\se)_{\alpha}\uua-(\se)_{\beta}\uub-N\un\Big)_{\beta\beta}\eqdef\cG_2.
\end{eqnarray*}
Consequently, we obtain
\begin{eqnarray}\label{eq:pressure}
&&-\Delta\Pi+4EH^2\Pi=H\Big(\Delta
H+\frac{H}{2E}\big((L-N)^2+4M^2\big)\Big)\nonumber\\
&&\qquad+\cG_1-2\varepsilon_0\cG_2-\frac{4\varepsilon_0
H}{E}\big(S_{11}L+2S_{12}M+S_{22}N\big)\eqdef\cG.
\end{eqnarray}

\begin{Remark}\label{rem:pressure}
In (\ref{eq:pressure}), $\mathcal{G}$ is a polynomial function of
$(\partial^kE,\partial^l\un,\partial^l\uua,\partial^l\uub,\partial^lL,\partial^lM,\partial^lN,\partial^lH)$,
where $0\leq |k|\leq3$ and $0\leq |l|\leq2$.
\end{Remark}

\begin{Remark}
It is reasonable that there is no term involving $\ww$ in $\cG$.
Actually, by differentiating the equation $\aaa^{\alpha}\cdot\vv_{,\alpha}=0$, and reformulating
the resulting equation in the isothermal coordinate, we can also
derive the equation of the pressure.
\end{Remark}

\setcounter{equation}{0}
\section{The linearized system}
In this section, we study the well-posedness of the linearized system of (\ref{eq:full system-R})-(\ref{eq:full system-velocity free}).
More precisely, we will consider the linear system
\begin{eqnarray}\label{eq:linearized system}
\left\{
\begin{array}{l}
\frac{\partial U_i}{\partial t}=\f 1 {E\sqrt{E}}\Delta(\sqrt{E}U_i)+G_i,\quad i=1,2,\\
\frac{\partial \un}{\partial t}=-\frac{1}{2E}\Delta H+B^1\cdot\na \un+G_3,\\
\frac{\partial{H}}{\partial{t}}=\frac{1}{2E}\Delta\un+B^2\cdot\na H+G_4,
\end{array}\right.
\end{eqnarray}
together with the initial condition
\ben\label{eq:initial condition-l}
(U_1,U_2,U^n,H)|_{t=0}=(U_1^0,U_2^0,U^n_0,H^0).
\een

Throughout this paper, we assume that $x=(x_1,x_2)\in \RT^2$.

\begin{Theorem}\label{thm:linear}
Let $s=2(k+1)$ for some integer $k\ge 2$, and let $T>0$. Suppose that
$E\in C([0,T];H^{s+1}(\RT^2))\cap C^1([0,T];H^{s-1}(\RT^2))$ and $E\ge c_0$
for some $c_0>0$. We also assume that $(U_1^0,U_2^0,U^n_0,H^0)\in H^{s-1}(\RT^2)$,
$(G_1,G_2)\in L^2(0,T;H^{s-2}(\RT^2)), (G_3,G_4)\in  L^2(0,T;H^{s-1}(\RT^2))$, and
$(B^1,B^2)\in L^2(0,T;H^{s-1}(\RT^2))$. Then there exists a unique solution $(U_1,U_2,U^n,H)$ on $[0,T]$ to
the linear system (\ref{eq:linearized system})-(\ref{eq:initial condition-l}) such that
\beno
&&(U_1,U_2)\in C([0,T];H^{s-1}(\RT^2))\cap L^2(0,T;H^s(\RT^2)),\\
&&(U^n,H)\in C([0,T];H^{s-1}(\RT^2)).
\eeno
Moreover, for any given $\veps>0$, it holds that
\ben\label{eq:energy}
&&E_s(t)\le C(\|E\|_{L^\infty_tH^{s-1}})\Big[E_s(0)+\int_0^t\cF_\veps(\|E\|_{H^{s+1}},\|E_t\|_{H^{s-1}})
(1+\|B\|_{H^s})E_s(\tau)d\tau\nonumber\\
&&\qquad\qquad+\int_0^t\|(G_1,G_2)(\tau)\|_{H^{s-2}}^2d\tau+\veps
\int_0^t\|(G_3,G_4)(\tau)\|_{H^{s-1}}^2d\tau\Big], \een
where $B=(B^1,B^2)$, $\cF_\veps$ is an increasing function, and $E_s(t)$ is defined by
\beno
E_s(t)\eqdef\|(U_1,U_2)(t)\|_{H^{s-1}}^2+\int_0^t\|(U_1,U_2)(\tau)\|_{H^s}^2d\tau+\|(U^n,H)(t)\|_{H^{s-1}}^2.
\eeno

\end{Theorem}

\no{\bf Proof.\,} The existence of $(U_1,U_2)$ is ensured by the classical parabolic theory,
whereas $(U^n,H)$  can be obtained by the duality method, see \cite{Ambrose} for example.
Here we only present the proof of the energy estimate. For this purpose,
let us introduce the energy functional $\cE$ defined by
\beno
\cE=\cE^1+\cE^2,
\eeno
with $\cE^1$ and $\cE^2$ given, respectively, by
\beno
&&\cE^1\eqdef\|\Lam^{s-1}U_{1}\|_{L^2}^2+\|\Lam^{s-1}U_{2}\|_{L^2}^2,\\
&&\cE^2\eqdef\|\Lambda(\f 1 E\Delta)^kU^n\|_{L^2}^2+\|\Lambda(\f 1 E\Delta)^kH\|_{L^2}^2,\quad \Lam=(-\Delta)^\f12.
\eeno

{\bf Step 1.}\textbf{ Estimate of $\cE^1$}

Taking the derivative to $\cE^1$ with respect to $t$, we obtain
\beno
\f12\f d {dt}\cE^1=\langle \Lam^{s-1}U_1, \Lam^{s-1}\p_tU_1\rangle+\langle \Lam^{s-1}U_2, \Lam^{s-1}\p_tU_2\rangle.
\eeno
By using the first equation of (\ref{eq:linearized system}), we get that
\beno
\langle \Lam^{s-1}U_1, \Lam^{s-1}\p_tU_1\rangle=
\langle \Lam^{s-1}U_1, \Lam^{s-1}\f 1 {E\sqrt{E}}\Delta(\sqrt{E}U_1)\rangle
+\langle \Lam^{s-1}U_1, \Lam^{s-1}G_1\rangle.
\eeno
By the Cauchy-Schwartz inequality, we have
\beno
\langle \Lam^{s-1}U_1, \Lam^{s-1}G_1\rangle\le \|U_1\|_{H^s}\|G_1\|_{H^{s-2}},
\eeno
and we write
\beno
&&\langle \Lam^{s-1}U_1, \Lam^{s-1}\f 1 {E\sqrt{E}}\Delta(\sqrt{E}U_1)\rangle\\
&&=\langle \Lam^{s-1}U_1, \big[\Lam^{s-1}, \f 1 {E\sqrt{E}}\Delta \sqrt{E}\big]U_1\rangle
+\langle \Lam^{s-1}U_1, \f 1 {E\sqrt{E}}\Delta(\sqrt{E}\Lam^{s-1}U_1)\rangle.
\eeno
By integration by parts and based on Lemmas \ref{lem:product} and \ref{lem:composition},
the second term of the right-hand side is bounded by
\beno
&&-c\|\Lam^{s}U_1\|^2_{L^2}+C\|\na E\|_{L^\infty}^2\|\Lam^{s-1}U_1\|^2_{L^2}
+C\|\na E\|_{L^\infty}\|\Lam^{s-1}U_1\|_{L^2}\|\Lam^{s}U_1\|_{L^2}\\
&&\le -\f c 2\|U_1\|^2_{H^{s}}+\cF(\|E\|_{H^{s+1}})\|U_1\|_{H^{s-1}}^2.
\eeno
Here the constant $c>0$ depends only on $c_0$.
And when Lemma \ref{lem:commutator} and Lemma \ref{lem:composition} are both used, the first term is bounded by
\beno
\cF(\|E\|_{H^{s+1}})\|U_1\|_{H^{s-1}}^2+\f c 4\|U_1\|_{H^{s}}^2,
\eeno
since we can write
\beno
\big[\Lam^{s-1}, \f 1 {E\sqrt{E}}\Delta \sqrt{E}\big]U_1=
\big[\Lam^{s-1}, \f 1 {E\sqrt{E}}\big]\Delta(\sqrt{E}U_1)+
\f 1 {E\sqrt{E}}\Delta\big[\Lam^{s-1},\sqrt{E}]U_1.
\eeno
Summing up the above estimates yields that
\beno
\langle \Lam^{s-1}U_1, \Lam^{s-1}\p_tU_1\rangle\le -c\|U_1\|^2_{H^s}
+\cF(\|E\|_{H^{s+1}})\|U_1\|_{H^{s-1}}^2+\|G_1\|_{H^{s-2}}^2,
\eeno
for some $c>0$. Similarly, we have
\beno
\langle \Lam^{s-1}U_2, \Lam^{s-1}\p_tU_2\rangle\le -c\|U_2\|^2_{H^s}
+\cF(\|E\|_{H^{s+1}})\|U_2\|_{H^{s-1}}^2+\|G_2\|_{H^{s-2}}^2.
\eeno
Hence, we obtain
\ben\label{eq:energy-1}
\f d {dt}\cE^1+c\|(U_1,U_2)\|^2_{H^s}\le \cF(\|E\|_{H^{s+1}})\|(U_1,U_2)\|_{H^{s-1}}^2
+\|(G_1,G_2)\|^2_{H^{s-2}}.
\een

{\bf Step 2.}\,\textbf{ Estimate of $\cE^2$}

Take the derivative to $\cE^2$ with respect to $t$ to obtain
\beno
&&\f12\f d {dt}\cE^2=\langle \Lam(\f 1 E\Delta)^k\p_tU^n, \Lam(\f 1 E\Delta)^kU^n\rangle
+\langle \Lam(\f 1 E\Delta)^k\p_tH, \Lam(\f 1 E\Delta)^kH\rangle\\
&&\quad+\langle\Lam\big[\p_t,(\f 1 E\Delta)^k\big]U^n, \Lam(\f 1 E\Delta)^kU^n\rangle
+\langle\Lam\big[\p_t,(\f 1 E\Delta)^k\big]H, \Lam(\f 1 E\Delta)^kH\rangle\\
&&\eqdef I+II+III+IV.
\eeno
Based on the last two equations of (\ref{eq:linearized system}), we get
\beno
&&I+II\\
&&=\langle \Lam(\f 1 E\Delta)^k(B^1\cdot\na U^n), \Lam(\f 1 E\Delta)^kU^n\rangle
+\langle \Lam(\f 1 E\Delta)^kG_3, \Lam(\f 1 E\Delta)^kU^n\rangle\\
&&\quad+\langle \Lam(\f 1 E\Delta)^k(B^2\cdot\na H), \Lam(\f 1 E\Delta)^kH\rangle
+\langle \Lam(\f 1 E\Delta)^kG_4,\Lam(\f 1 E\Delta)^kH\rangle\\
&&\eqdef I_1+I_2+II_1+II_2.
\eeno
Here we use the following fact:
\beno
\langle \Lam(\f 1 E\Delta)^k(-\f 1 {2E}\Delta)H, \Lam(\f 1 E\Delta)^kU^n\rangle
+\langle \Lam(\f 1 E\Delta)^k(\f 1 {2E}\Delta)U^n, \Lam(\f 1 E\Delta)^kH\rangle=0.
\eeno
Using Lemma \ref{lem:operator-commuatator-time} and Lemma \ref{lem:operator-upper bound}, we get
\beno
|III|+|IV|\le \cF(\|(E,E_t)\|_{H^{s-1}})\big(\|U^n\|_{H^{s-1}}^2+\|H\|_{H^{s-1}}^2\big),
\eeno
and by Lemma \ref{lem:operator-upper bound},
\beno
|I_2|+|II_2|\le \cF_\veps(\|E\|_{H^{s-1}})\big(\|U^n\|_{H^{s-1}}^2+\|H\|_{H^{s-1}}^2\big)
+\veps\|(G_3,G_4)\|_{H^{s-1}}^2.
\eeno
To estimate $I_1$, we write
\beno
I_1&=&\langle\Lam(\f 1 E\Delta)^k(B^1\cdot\na U^n)-B^1\cdot\na \Lam(\f 1 E\Delta)^kU^n, \Lam(\f 1 E\Delta)^kU^n\rangle\\
&&+\langle B^1\cdot\na \Lam(\f 1 E\Delta)^kU^n, \Lam(\f 1 E\Delta)^kU^n\rangle
=I_{11}+I_{12}.
\eeno
We have by Lemma \ref{lem:operator-upper bound} that
\beno
|I_{12}|\le \cF(\|E\|_{H^{s-1}})\|\na B^1\|_{L^\infty}\|U^n\|_{H^{s-1}}^2.
\eeno
We further write
\beno
I_{11}&=&\langle\Lam\big[(\f 1 E\Delta)^k,B^1\big]\cdot\na U^n
+\big[\Lam, B^1\big]\cdot\na(\f 1 E\Delta)^kU^n, \Lam(\f 1 E\Delta)^kU^n\rangle\\
&&+\langle B^1\cdot\Lam\big[(\f 1 E\Delta)^k,\na \big]U^n, \Lam(\f 1 E\Delta)^kU^n\rangle,
\eeno
which along with Lemma \ref{lem:operator-upper bound} and
Lemma \ref{lem:operator-commuatator-time}-\ref{lem:operator-commuatator-function} implies that
\beno
|I_{11}|\le \cF(\|E\|_{H^{s}})\|B^1\|_{H^{s-1}}\|U^n\|_{H^{s-1}}^2.
\eeno

On the basis of the above estimates, we obtain
\ben\label{eq:energy-2}
\f12\f d {dt}\cE^2&\le& \cF_\veps(\|E\|_{H^{s}},\|E_t\|_{H^{s-1}})(1+\|B\|_{H^{s-1}})\big(\|U^n\|_{H^{s-1}}^2+\|H\|_{H^{s-1}}^2\big)\nonumber\\
&&\quad+\veps\|(G_3,G_4)\|_{H^{s-1}}^2.
\een

{\bf Step 3. $L^2$ estimate}

Taking the $L^2$ energy estimate for $U_i(i=1,2)$, we obtain
\beno
&&\f 12 \f d {dt}\|U_i\|_{L^2}^2+\|\f 1 E\na(\sqrt{E}U_i)\|_{L^2}^2\\
&&\le C(\|E\|_{L^\infty})\|\f 1 E\na E\|_{L^\infty}\|\na(\sqrt{E}U_i)\|_{L^2}\|U_i\|_{L^2}
+\|G_i\|_{L^2}\|U_i\|_{L^2}\\
&&\le \cF(\|E\|_{H^3})\|U_i\|_{L^2}^2+\|G_i\|_{L^2}^2+\|\f 1 E\na(\sqrt{E}U_i)\|_{L^2}^2.
\eeno
Taking the $L^2$ energy estimate for $(U^n,H)$, we get
\beno
&&\langle E\p_tU^n,U^n\rangle+\langle E\p_tH,H\rangle\\
&&=\langle EB^1\cdot\na U^n,U^n\rangle+\langle EG_3,U^n\rangle
+\langle EB^2\cdot\na H,H\rangle+\langle EG_4,H\rangle,
\eeno
from which, we infer that
\beno
&&\f12\f d {dt}\big(\|\sqrt{E}U^n\|_{L^2}^2+\|\sqrt{E}H\|_{L^2}^2\big)\\
&&\le \cF_\veps(\|(E,E_t)\|_{H^3})(1+\|B\|_{H^3})\big(\|U^n\|_{L^2}^2+\|H\|_{L^2}^2\big)
+\veps\|(G_3,G_4)\|_{L^2}^2.
\eeno
Thus, we obtain
\ben\label{eq:energy-L2}
\f d {dt}\cE_0&\le& \cF(\|E\|_{H^3})\|(U_1,U_2)\|_{L^2}^2
+\cF_\veps(\|(E,E_t)\|_{H^3})(1+\|B\|_{H^3})\|(U^n,H)\|_{L^2}^2\nonumber\\
&&+\|(G_1,G_2)\|_{L^2}^2+\veps\|(G_3,G_4)\|_{L^2}^2.
\een
Here $\cE_0\eqdef \|(U_1,U_2)\|_{L^2}^2+\|\sqrt{E}U^n\|_{L^2}^2+\|\sqrt{E}H\|_{L^2}^2.$
\vspace{0.1cm}

Now we are in position to complete the proof.
Taken together (\ref{eq:energy-1})-(\ref{eq:energy-L2}) yields that
\beno
\f d {dt}(\cE+\eta\cE_0)+c\|(U_1,U_2)\|^2_{H^s}&\le& \cF_\veps(\|E\|_{H^{s+1}},\|E_t\|_{H^{s-1}})(1+\|B\|_{H^s})(\cE+\eta\cE_0)\\
&&+\eta\big(\|(G_1,G_2)\|_{H^{s-2}}^2+\veps\|(G_3,G_4)\|_{H^{s-1}}^2\big),
\eeno which implies (\ref{eq:energy}) by taking $\eta$ to be bigger than
$C(\|E\|_{L^\infty_tH^{s-1}})$, since we have by Lemma \ref{lem:operator-lowerbound}
and an interpolation argument that
\beno
\cE^2\ge c\|(U^n,H)\|_{H^{s-1}}^2-C(\|E\|_{H^{s-1}})\|(U^n,H)\|_{L^2}^2.
\eeno
This completes the proof of Theorem \ref{thm:linear}. \endproof

\setcounter{equation}{0}
\section{Nonlinear system}

This section is devoted to solving the nonlinear system (\ref{eq:full system-R})-(\ref{eq:full system-velocity free}).

\subsection{Iteration scheme}
We will construct the solution $(\RR, U_1,U_2,U^n,H)$ by the iteration method. First of all,
we take
\beno
(\RR^{(0)},U_1^{(0)}, U_2^{(0)}, U^{n, (0)},E^{(0)}, H^{(0)})=(\RR_0, U_1^0,U_2^0,U^n_0,E_0, H_0);
\eeno
And, $(W_1^{(0)},W^{(0)}_2)$ are determined by solving the following elliptic system:
\beno
\left\{
\begin{array}{l}
\Big(\frac{W_1^{(0)}}{\sqrt{E^{(0)}}}\Big)_{\alpha}-\Big(\frac{W_2^{(0)}}{\sqrt{E^{(0)}}}\Big)_{\beta}=\frac{U^{n,(0)}(L^{(0)}-N^{(0)})}{E^{(0)}},\\
\Big(\frac{W_1^{(0)}}{\sqrt{E^{(0)}}}\Big)_{\beta}+\Big(\frac{W_2^{(0)}}{\sqrt{E^{(0)}}}\Big)_{\alpha}=\frac{2U^{n,(0)} M^{(0)}}{E^{(0)}}.
\end{array}\right.
\eeno
The pressure $\Pi^{(0)}$ is given by
\beno
-\Delta\Pi^{(0)}+4E^{(0)}(H^{(0)})^2\Pi^{(0)}=\cG^{(0)}
\eeno
with $\cG^{(0)}$ determined by $(\RR^{(0)},U_1^{(0)}, U_2^{(0)}, U^{n, (0)})$ \big(see (\ref{eq:pressure})\big).

Assume that $(U_1^{(\ell)}, U_2^{(\ell)}, U^{n, (\ell)}, H^{(\ell)},E^{(\ell)}, W_1^{(\ell)}, W_2^{(\ell)},\RR^{(\ell)})$ has been constructed.
We denote
\beno
L^{(\ell)}=\RR_{\al\al}^{(\ell)}\cdot\nn^{(\ell)},\quad M^{(\ell)}=\RR_{\al\be}^{(\ell)}\cdot\nn^{(\ell)},\quad
N^{(\ell)}=\RR_{\be\be}^{(\ell)}\cdot\nn^{(\ell)},
\quad\nn^{(\ell)}=\frac{\RR_{\alpha}^{(\ell)}\times\RR_{\beta}^{(\ell)}}{|\RR_{\alpha}^{(\ell)}\times\RR_{\beta}^{(\ell)}|}.
\eeno
Then we construct $(U_1^{(\ell+1)}, U_2^{(\ell+1)}, U^{n, (\ell+1)}, H^{(\ell+1)})$
by solving the following linear system:
\begin{eqnarray}\label{eq:approximate system-velocity}
\left\{
\begin{array}{l}
\frac{\partial \uua^{(\ell+1)}}{\partial t}=\epsilon_0\f 1 {E^{(\ell)}\sqrt{E^{(\ell)}}}\Delta(\sqrt{E^{(\ell)}}\uua^{(\ell+1)})+\cF_1^{(\ell)},\\
\frac{\partial \uub^{(\ell+1)}}{\partial t}=\epsilon_0\f 1 {E^{(\ell)}\sqrt{E^{(\ell)}}}\Delta(\sqrt{E^{(\ell)}}\uub^{(\ell+1)})+\cF_2^{(\ell)},\\
\frac{\partial U^{n,(\ell+1)}}{\partial t}=-\frac{1}{2E^{(\ell)}}\Delta
H^{(\ell+1)}-(2u_i^{(\ell)}-w_i^{(\ell)})\p_iU^{n,(\ell+1)}+\cF_3^{(\ell)},\\
\frac{\partial{H^{(\ell+1)}}}{\partial{t}}=\frac{1}{2E^{(\ell)}}\Delta U^{n,(\ell+1)}
+w_i^{(\ell)}\p_iH_{\alpha}^{(\ell+1)}+F_7^{(\ell)},\\
(U_1^{(\ell+1)},U_2^{(\ell+1)},U^{n,(\ell+1)},H^{(\ell+1)})|_{t=0}=(U_1^0,U_2^0,U^n_0,H_0),
\end{array}\right.
\end{eqnarray}
where  $u_i^{(\ell)}=\frac{U_i^{(\ell)}}{\sqrt{E^{(\ell)}}},\, w_i^{(\ell)}=\frac{W_i^{(\ell)}}{\sqrt{E^{(\ell)}}}$ for $i=1,2$ and
\beno
\cF^{(\ell)}_1&=&-\f 1 {\sqrt{E^{(\ell)}}}{\Pi_{\alpha}^{(\ell)}}
+\f {2\epsilon_0}{E^{(\ell)}\sqrt{E^{(\ell)}}}F_{1\alpha}^{(\ell)}+\frac 1{2\sqrt{E^{(\ell)}}}F_{2\alpha}^{(\ell)}\\
&&-\f {\epsilon_0}{E^{(\ell)}\sqrt{E^{(\ell)}}}F_{3\beta}^{(\ell)}
-(u_1^{(\ell)}-w_1^{(\ell)})\uaa^{(\ell)}-(u_2^{(\ell)}-w_2^{(\ell)})\uab^{(\ell)}+F_4^{(\ell)},\\
\cF^{(\ell)}_2&=&-\f 1 {\sqrt{E^{(\ell)}}}{\Pi_{\be}^{(\ell)}}
+\f {2\epsilon_0}{E^{(\ell)}\sqrt{E^{(\ell)}}}F_{1\be}^{(\ell)}+\frac 1{2\sqrt{E^{(\ell)}}}F_{2\be}^{(\ell)}\\
&&-\f {\epsilon_0}{E^{(\ell)}\sqrt{E^{(\ell)}}}F_{3\beta}^{(\ell)}
-(u_1^{(\ell)}-w_1^{(\ell)})\uba^{(\ell)}-(u_2^{(\ell)}-w_2^{(\ell)})\ubb^{(\ell)}+F_5^{(\ell)},\\
\cF_3^{(\ell)}&=& 2H^{(\ell)}\Pi^{(\ell)}
+\f {2\epsilon_0}{E^{(\ell)}\sqrt{E^{(\ell)}}}(L^{(\ell)}\uaa^{(\ell)}+M^{(\ell)}\uab^{(\ell)}+M^{(\ell)}\uba^{(\ell)}+N^{(\ell)}\ubb^{(\ell)})+F_6^{(\ell)},
\eeno
with $F_i^{(\ell)}(i=1,\cdots,7)$ are given in Section 3.3 where $(U_1,U_2,U^n,W_1,W_2,E,L,M,N)$
are replaced by $(U_1^{(\ell)},U_2^{(\ell)},U^{n,(\ell)},W_1^{(\ell)},W_2^{(\ell)},E^{(\ell)},L^{(\ell)},M^{(\ell)},N^{(\ell)})$.
Let $\widetilde{\RR}_t^{{k+1}}$ be given by
\ben\label{eq:R-app-tilde}
\widetilde{\RR}_t^{(\ell+1)}=U^{n,(\ell)}\frac{\RR_{\alpha}^{(\ell)}\times \RR_\be^{(\ell)}}{E^{(\ell)}}+W_1^{(\ell)}\frac{\RR_{\alpha}^{(\ell)}}{\sqrt{E^{(\ell)}}}
+W_2^{(\ell)}\frac{\RR_{\be}^{(\ell)}}{\sqrt{E^{(\ell)}}},\quad \widetilde{\RR}|_{t=0}=\RR_0.
\een
And, $\widehat{\RR}^{(\ell+1)}$ is determined by solving
\ben\label{eq:R-hat}
\Delta \widehat{\RR}^{(\ell+1)}-\widehat{\RR}^{(\ell+1)}
=2H^{(\ell)}\widetilde{\RR}_\al^{(\ell+1)}\times\widetilde{\RR}_\be^{(\ell+1)}-\widetilde{\RR}^{(\ell+1)}.
\een
Then we construct the surface $\RR^{(\ell+1)}$ by solving the following elliptic equation:
\ben\label{eq:R-app}
\Delta \RR^{(\ell+1)}-\RR^{(\ell+1)}=2H^{(\ell)}\widehat{\RR}_\al^{(\ell+1)}\times\widehat{\RR}_\be^{(\ell+1)}
-\widetilde{\RR}^{(\ell+1)}.
\een

Next we define $E^{(\ell+1)}$ by solving
\begin{eqnarray}\label{eq:E-app}
\Delta E^{(\ell+1)}-2E^{(\ell+1)}&=&2(\RR_{\alpha\beta}^{(\ell+1)}\cdot\RR_{\alpha\beta}^{(\ell+1)}
-\RR_{\alpha\alpha}^{(\ell+1)}\cdot\RR_{\beta\beta}^{(\ell+1)})\nonumber\\
&&\quad-\big(\RR_\al^{(\ell+1)}\cdot\RR_\al^{(\ell+1)}+\RR_\be^{(\ell+1)}\cdot\RR_\be^{(\ell+1)}\big).
\end{eqnarray}
And, $(W_1^{(\ell+1)},W^{(\ell+1)}_2)$ is determined by solving
\ben\label{eq:W-app}
\left\{
\begin{array}{l}
\Big(\frac{W_1^{(\ell+1)}}{\sqrt{E^{(\ell)}}}\Big)_{\alpha}-\Big(\frac{W_2^{(\ell+1)}}{\sqrt{E^{(\ell)}}}\Big)_{\beta}
=\frac{U^{n,(\ell)}(L^{(\ell)}-N^{(\ell)})}{E^{(\ell)}},\\
\Big(\frac{W_1^{(\ell+1)}}{\sqrt{E^{(\ell)}}}\Big)_{\beta}+\Big(\frac{W_2^{(\ell+1)}}{\sqrt{E^{(\ell)}}}\Big)_{\alpha}=\frac{2U^{n,(\ell)} M^{(\ell)}}{E^{(\ell)}}.
\end{array}\right.
\een
Finally, we define the pressure $\Pi^{(\ell+1)}$ by solving
\begin{eqnarray}\label{eq:Pressure-app}
-\Delta\Pi^{(\ell+1)}+4E^{(\ell)}(H^{(\ell)})^2\Pi^{(\ell+1)}=\cG^{(\ell)},
\end{eqnarray}
with $\cG^{(\ell)}$ determined by $(U_1^{(\ell)}, U_2^{(\ell)}, U^{n, (\ell)},E^{(\ell)},L^{(\ell)},M^{(\ell)},N^{(\ell)})$, see (\ref{eq:pressure}).

\begin{Remark} If $\RR^{(l+1)}$ is directly defined by (\ref{eq:R-app-tilde}),
then we can only obtain the $H^{s-1}$ regularity of $\RR^{(l+1)}$.
However, we need the $H^{s+1}$ regularity of $\RR^{(l+1)}$ to close the energy estimates. Motivated by (\ref{eq:R-elliptic equation}), we determine
$\RR^{(l+1)}$ by using (\ref{eq:R-app-tilde})-(\ref{eq:R-app}) so that
the $H^{s+1}$ regularity of $\RR^{(l+1)}$ can be obtained by the elliptic estimates.
\end{Remark}

\subsection{Nonlinear estimates}

Before presenting the estimates, let us make the following assumptions on the step-$\ell^{th}$
approximate solutions $(\RR^{(\ell)},U_1^{(\ell)}, U_2^{(\ell)}, U^{n, (\ell)},E^{(\ell)},H^{(\ell)})$:
\ben
&&E^{(\ell)}(x,t)\ge c_0>0,  \quad \textrm{for any } (t,x)\in [0,T]\times \RT^2,\label{eq:E-ass}\\
&&\int_{\RT^2}(H^{(\ell)}(x,t))^2dx\ge c_1>0,\quad \textrm{for any } t\in [0,T],\label{eq:H-ass}\\
&&|\RR_{\al}^{(\ell)}(x,t)\times\RR_{\be}^{(\ell)}(x,t)|\ge c_0,\quad \textrm{for any } (t,x)\in [0,T]\times \RT^2\label{eq:n-ass}\\
&&\|(U_1^{(\ell)},U_2^{(\ell)}\|_{L^\infty(0,T;H^{s-1})}+\|(U_1^{(\ell)},U_2^{(\ell)})\|_{L^2(0,T;H^{s})}\le C_1,\label{eq:U-ass}\\
&&\|(U^{n,(\ell)},H^{(\ell)})\|_{L^\infty(0,T;H^{s-1})}\le C_2,\label{eq:Un-ass}\\
&&\|\RR^{(\ell)}\|_{C^i([0,T];H^{s+1-2i})}+\|E^{(\ell)}\|_{C^i([0,T];H^{s+1-2i})}\le C_3,\quad i=0,1,\label{eq:Rh-ass}\\
&&\|\RR^{(\ell)}\|_{L^\infty(0,T;H^{s-1})}+\|E^{(\ell)}\|_{L^\infty(0,T;H^{s-1})}\le C_4.\label{eq:R-ass}
\een
Here $T>0, s=2(k+1), k\ge 2$, and $C_1,C_2,C_3$, and $C_4$ are some fixed constants to be determined in Section 5.3.
Note that the assumptions (\ref{eq:E-ass})-(\ref{eq:n-ass}) and (\ref{eq:Rh-ass})-(\ref{eq:R-ass}) are made so that
we can use Theorem \ref{thm:linear} at each step of the iterations,
and the assumptions (\ref{eq:U-ass})-(\ref{eq:Un-ass}) are determined by
the energy estimates for the linearized system.

In what follows, we denote $\cC$ by an increasing function, which may be different from line to line.
From the definition, it is easy to see that
\ben\label{eq:second fundamental-est}
\|(L^{(\ell)},M^{(\ell)},N^{(\ell)})\|_{L^\infty(0,T;H^{s-1})}\le \cC(C_3).
\een
Using Lemma \ref{lem:elliptic system} and Lemma \ref{lem:product}, we find that
\ben\label{eq:W-est}
\|W_1^{(\ell+1)}\|_{L^\infty(0,T;H^{s})}+\|W_2^{(\ell+1)}\|_{L^\infty(0,T;H^{s})}\le \cC(C_2,C_3).
\een
From (\ref{eq:second fundamental-est}), (\ref{eq:W-est}), and Lemmas \ref{lem:product}-\ref{lem:composition}, we infer that
for $i=1,\cdots,7$,
\ben\label{eq:F-est}
\|F_i^{(\ell)}\|_{L^\infty(0,T;H^{s-1})}\le \cC(C_1,C_2,C_3).
\een
Thanks to Remark \ref{rem:pressure}, we get by using Lemma \ref{lem:product}-\ref{lem:composition} that
\beno
\|\cG^{(\ell)}\|_{L^\infty(0,T;H^{s-3})}\le \cC(C_1,C_2,C_3).
\eeno
Thus, we infer from Lemma \ref{lem:elliptic equation} that
\ben\label{eq:pressure-est}
\|\Pi^{(\ell+1)}\|_{L^\infty(0,T;H^{s-1})}\le \cC(C_1,C_2,C_3).
\een

Using (\ref{eq:second fundamental-est})-(\ref{eq:pressure-est}), we obtain
\begin{Proposition}\label{prop:nonlinear estimate}
The nonlinear terms $\cF^{(\ell)}_1, \cF^{(\ell)}_2$, and $\cF^{(\ell)}_3$ satisfy
\beno
&&\|\cF^{(\ell)}_i\|_{L^\infty(0,T;H^{s-2})}\le \cC(C_1,C_2,C_3),\quad i=1,2,\\
&&\|\cF^{(\ell)}_3\|_{L^2(0,T;H^{s-1})}\le \cC(C_1,C_2,C_3).
\eeno
\end{Proposition}

In order to prove the convergence of the iteration scheme, we need to establish some difference estimates
in the lower-order Sobolev spaces. For this, we set
\beno
&&\delta^\ell_{U_i}=U_i^{(\ell)}-U_i^{(\ell-1)}(i=1,2),\quad \delta^\ell_{U^n}=U^{n,(\ell)}-U^{n,(\ell-1)},
\quad \delta^\ell_{H}=H^{(\ell)}-H^{(\ell-1)},\\
&&\delta_\RR^\ell=\RR^{(\ell)}-\RR^{(\ell-1)},\quad\delta_E^\ell=E^{(\ell)}-E^{(\ell-1)}.
\eeno

First of all, we have
\beno
&&\|(L^{(\ell)},M^{(\ell)},N^{(\ell)})-(L^{(\ell-1)},M^{(\ell-1)},N^{(\ell-1)})\|_{H^{s-3}}
\le \cC(C_3)\|\delta_\RR^\ell\|_{H^{s-1}},\\
&&\|W^{(\ell+1)}-W^{(\ell)}\|_{H^{s-2}}\le \cC(C_3)(\|\delta_\RR^\ell\|_{H^{s-1}}+\|\delta_E^\ell\|_{H^{s-1}}+\|\delta_{U^n}^\ell\|_{H^{s-3}}),
\eeno
which imply that for $i=1,\cdots,7$,
\beno
\|F_i^{(\ell)}-F_i^{(\ell-1)}\|_{H^{s-3}}&\le& \cC(C_1,C_2,C_3)\big(\|\delta_\RR^\ell\|_{H^{s-1}}+\|\delta_E^\ell\|_{H^{s-1}}\\
&&+\|\delta^\ell_{U_1}\|_{H^{s-3}}+\|\delta^\ell_{U_2}\|_{H^{s-3}}+\|\delta^\ell_{U^n}\|_{H^{s-3}}+\|\delta^\ell_{H}\|_{H^{s-3}}\big).
\eeno
Similarly, we can obtain
\beno
\|\cG^{(\ell)}-\cG^{(\ell-1)}\|_{H^{s-5}}&\le& \cC(C_1,C_2,C_3)\big(\|\delta_\RR^\ell\|_{H^{s-1}}+\|\delta_E^\ell\|_{H^{s-1}}\\
&&+\|\delta^\ell_{U_1}\|_{H^{s-3}}+\|\delta^\ell_{U_2}\|_{H^{s-3}}+\|\delta^\ell_{U^n}\|_{H^{s-3}}+\|\delta^\ell_{H}\|_{H^{s-3}}\big).
\eeno
Hence, we infer from Lemma \ref{lem:elliptic equation} that
\beno
\|\Pi^{(\ell+1)}-\Pi^{(\ell)}\|_{H^{s-3}}&\le& \cC(C_1,C_2,C_3)\big(\|\delta_\RR^\ell\|_{H^{s-1}}+\|\delta_E^\ell\|_{H^{s-1}}\\
&&+\|\delta^\ell_{U_1}\|_{H^{s-3}}+\|\delta^\ell_{U_2}\|_{H^{s-3}}+\|\delta^\ell_{U^n}\|_{H^{s-3}}+\|\delta^\ell_{H}\|_{H^{s-3}}\big).
\eeno

From the above estimates, we can deduce

\begin{Proposition}\label{prop:nonlinear estimate-difference}
For $i=1,2$, it holds that
\beno
\|\cF_i^{(\ell)}-\cF_i^{(\ell-1)}\|_{H^{s-4}}&\le& \cC(C_1,C_2,C_3)\big(\|\delta_\RR^\ell\|_{H^{s-1}}+\|\delta_E^\ell\|_{H^{s-1}}\\
&&+\|\delta^\ell_{U_1}\|_{H^{s-3}}+\|\delta^\ell_{U_2}\|_{H^{s-3}}+\|\delta^\ell_{U^n}\|_{H^{s-3}}+\|\delta^\ell_{H}\|_{H^{s-3}}\big),\\
\|\cF_3^{(\ell)}-\cF_3^{(\ell-1)}\|_{H^{s-3}}&\le& \cC(C_1,C_2,C_3)\big(\|\delta_\RR^\ell\|_{H^{s-1}}+\|\delta_E^\ell\|_{H^{s-1}}\\
&&+\|\delta^\ell_{U_1}\|_{H^{s-2}}+\|\delta^\ell_{U_2}\|_{H^{s-2}}+\|\delta^\ell_{U^n}\|_{H^{s-3}}+\|\delta^\ell_{H}\|_{H^{s-3}}\big).
\eeno
\end{Proposition}

\subsection{Proof of the main result}

To simplify the analysis, we will first prove the well-posedness of the system by assuming that
the surface can be globally parameterized by the isothermal coordinates.
In Section 5.5, we will indicate how to extend this result to a general closed
surface, and thereby conclude the proof of Theorem \ref{thm:main}.

\begin{Theorem}
Let $s=2(k+1)$ for some integer $k\ge 2$. Assume that $(U_1^0,U_2^0,U^n_0)\in H^{s-1}(\RT^2)$,
and the initial surface $\RR_0\in H^{s+1}$. Moreover, the coefficient of the first fundamental form $E_0$ and
the mean curvature $H_0$ satisfy
\beno
&&(\sqrt{E_0}\uua^0)_{\alpha}+(\sqrt{E_0}\uub^0)_{\beta}-2E_0H_0\un_0=0,\\
&&E_0(x)\ge 2c_0,\quad \int_{\RT^2} H_0^2(x)dx\ge 2c_1,
\eeno
for some $c_0>0, c_1>0$. Then there exists $T>0$ such that the nonlinear
system (\ref{eq:full system-R})-(\ref{eq:full system-velocity free})
has a unique solution $(\RR, U_1,U_2,U^n, H)$ on $[0,T]$ satisfying
\beno
&&(U_1,U_2)\in C([0,T];H^{s-1})\cap L^2(0,T;H^s),\\
&&\RR\in C([0,T];H^{s+1}), \quad (U^n,H) \in C([0,T];H^{s-1}).
\eeno
\end{Theorem}

\begin{Remark}
We have chosen the isothermal coordinate for the initial surface. Hence, the conditions
\beno
&E_0(x)\ge 2c_0,\quad \int_{\RT^2} H_0^2(x)dx\ge 2c_1 \eeno are
naturally satisfied for any smooth closed surface.
\end{Remark}

\no{\bf Proof.\,} We  split the proof into two steps.\vspace{0.1cm}

{\bf Step 1.\,} \textbf{Uniform estimates}

Let us assume that $(\RR^{(\ell)},U_1^{(\ell)}, U_2^{(\ell)}, U^{n, (\ell)},E^{(\ell)},H^{(\ell)})$ satisfies
(\ref{eq:E-ass})-(\ref{eq:R-ass}). We will show
that $(\RR^{(\ell+1)},U_1^{(\ell+1)}, U_2^{(\ell+1)}, U^{n, (\ell+1)},E^{(\ell+1)},H^{(\ell+1)})$
also satisfies the same estimates.

We denote
\beno
E^{(\ell)}_s(t)\eqdef \|(U_1^{(\ell)},U_2^{(\ell)})(t)\|_{H^{s-1}}^2
+\int_0^t\|(U_1^{(\ell)},U_2^{(\ell)})(\tau)\|_{H^s}^2d\tau+\|(U^{n,(\ell)},H^{(\ell)})(t)\|_{H^{s-1}}^2.
\eeno
Then we infer from Theorem \ref{thm:linear} and Proposition \ref{prop:nonlinear estimate} that
\beno
E^{(\ell+1)}_s(t)&\le& \cC(C_4)\Big(E_0+\cC_\veps(C_1,C_2,C_3)\int_0^t(1+\|(U_1^{(\ell)},U_2^{\ell})\|_{H^s})E^{(\ell+1)}_s(\tau)d\tau\\
&&+\cC(C_1,C_2,C_3)(t+\veps)\Big).
\eeno
Here $E_0\eqdef \|(U_1^0,U_2^0,U^n_0,H^0)\|_{H^{s-1}}$. Then we get by Gronwall's inequality that
\beno
E^{(\ell+1)}_s(t)\le \big(\cC(C_4)E_0+\cC(C_1,C_2,C_3,C_4)(t+\veps)\big)\exp(\cC(C_1,C_2,C_3,C_4)t).
\eeno
Taking $T$ and $\veps$ small enough yields that
\beno
E^{(\ell+1)}_s(t)\le 2\cC(C_4)E_0\quad\textrm{ for }t\in [0,T].
\eeno
This means that if we take $C_1=C_2=2\cC(C_4)E_0$, $(U_1^{(\ell+1)}, U_2^{(\ell+1)}, U^{n, (\ell+1)}, H^{(\ell+1)}))$
satisfies (\ref{eq:U-ass})-(\ref{eq:Un-ass}).

Due to (\ref{eq:R-app-tilde}), we find that
\beno
&&\|\widetilde{\RR}^{(\ell+1)}\|_{L^\infty(0,T;H^{s-1})}\le \|R_0\|_{H^{s-1}}
+\cC(C_1,C_2,C_3)t.
\eeno
Hence, by taking $T$ to be small enough if necessary, we get
\ben\label{eq:RE-est}
\|\widetilde{\RR}^{(\ell+1)}\|_{L^\infty(0,T;H^{s-1})}\le 2\|\RR_0\|_{H^{s-1}}.
\een
We also have by (\ref{eq:R-app-tilde}) that
\ben\label{eq:RE-est-time}
\|\p_t\widetilde{\RR}^{(\ell+1)}\|_{L^\infty(0,T;H^{s-2})}\le \cC(C_2,C_4).
\een
We get by the elliptic estimate that
\beno
\|\widehat{\RR}^{(\ell+1)}\|_{L^\infty(0,T;H^{s})}\le \cC(\|\RR_0\|_{H^{s-1}}, C_2),
\eeno
which along with (\ref{eq:RE-est}) and (\ref{eq:R-app}) implies that
\beno
\|\RR^{(\ell+1)}\|_{L^\infty(0,T;H^{s+1})}\le \cC(\|\RR_0\|_{H^{s}},C_2).
\eeno
Hence, by (\ref{eq:E-app}) and the elliptic estimate,
\beno
\|E^{(\ell+1)}\|_{L^\infty(0,T;H^{s+1})}\le \cC(\|\RR_0\|_{H^{s}},C_2).
\eeno
Taking the derivative to (\ref{eq:R-app}) and (\ref{eq:E-app}) with respect to time,
we get by (\ref{eq:RE-est-time}) that
\ben\label{eq:R-est-time}
\|\p_t{\RR}^{(\ell+1)}\|_{L^\infty(0,T;H^{s-1})}
+\|\p_t{E}^{(\ell+1)}\|_{L^\infty(0,T;H^{s-1})}\le \cC(\|\RR_0\|_{H^{s}},C_2,C_4).
\een
Hence, taking $C_3=\cC(\|\RR_0\|_{H^{s}},C_2,C_4)$, we see that $(\RR^{(\ell+1)},E^{(\ell+1)})$
satisfies (\ref{eq:Rh-ass}).

Now taking $C_4=2(\|\RR_0\|_{H^{s-1}}+\|E_0\|_{H^{s-1}})\le C\|\RR_0\|_{H^s}$, it follows from (\ref{eq:R-est-time}) that
\beno
\|({\RR}^{(\ell+1)},{E}^{(\ell+1)})\|_{L^\infty(0,T;H^{s-1})}
\le \|\RR_0\|_{H^{s-1}}+\|E_0\|_{H^{s-1}}
+\cC(\|\RR_0\|_{H^{s}},C_2,C_4)T,
\eeno
which implies that $(\RR^{(\ell+1)},E^{(\ell+1)})$ satisfies (\ref{eq:R-ass}) when $T$ is taken to be small enough.
Similarly, we can show that $(\RR^{(\ell+1)},E^{(\ell+1)},H^{(\ell+1)})$ also satisfies (\ref{eq:E-ass})-(\ref{eq:n-ass}).

In conclusion, we prove that there exists a $T>0$ depending only on $\|(U_1^0,U_2^0,U^n_0)\|_{H^{s-1}}$
and $\|\RR_0\|_{H^{s+1}}$ such that  (\ref{eq:E-ass})-(\ref{eq:R-ass})
hold for $(\RR^{(\ell+1)},U_1^{(\ell+1)}, U_2^{(\ell+1)}, U^{n, (\ell+1)},E^{(\ell+1)},H^{(\ell+1)})$.
\vspace{0.1cm}

{\bf Step 2.\,} \textbf{Existence and uniqueness}

It suffices to show that the approximate solution sequence is a Cauchy sequence.
For this purpose, we set
\beno
&&\delta^{\ell+1}_{U_i}=U_i^{(\ell+1)}-U_i^{(\ell)}(i=1,2),\quad \delta^{\ell+1}_{U^n}=U^{n,(\ell+1)}-U^{n,(\ell)},
\quad \delta^{\ell+1}_{H}=H^{(\ell+1)}-H^{(\ell)},\\
&&\delta_\RR^{\ell+1}=\RR^{(\ell+1)}-\RR^{(\ell)},\quad\delta_E^{\ell+1}=E^{(\ell+1)}-E^{(\ell)}.
\eeno
Then $(\delta^{\ell+1}_{U_1},\delta^{\ell+1}_{U_2},\delta^{\ell+1}_{U_n},\delta^{\ell+1}_{H})$ satisfies the following system:
\begin{eqnarray}\label{eq:difference system-velocity}
\left\{
\begin{array}{l}
\frac{\partial\delta^{\ell+1}_{U_1}}{\partial t}=\epsilon_0\f 1 {E^{(\ell)}\sqrt{E^{(\ell)}}}\Delta(\sqrt{E^{(\ell)}}\delta^{\ell+1}_{U_i})+\delta\cF_1^{(\ell)},\\
\frac{\partial \delta^{\ell+1}_{U_2}}{\partial t}=\epsilon_0\f 1 {E^{(\ell)}\sqrt{E^{(\ell)}}}\Delta(\sqrt{E^{(\ell)}}\delta^{\ell+1}_{U_2})+\delta\cF_2^{(\ell)},\\
\frac{\partial \delta^{\ell+1}_{U^n}}{\partial t}=-\frac{1}{2E^{(\ell)}}\Delta
\delta^{\ell+1}_{H}-(2u_i^{(\ell)}-w_i^{(\ell)})\p_i\delta^{\ell+1}_{U^n}+\delta\cF_3^{(\ell)},\\
\frac{\partial{\delta^{\ell+1}_{H}}}{\partial{t}}=\frac{1}{2E^{(\ell)}}\Delta \delta^{\ell+1}_{U^n}
+w_i^{(\ell)}\p_i\delta^{\ell+1}_{H}+\delta\cF_4^{(\ell)},\\
(\delta^{\ell+1}_{U_1},\delta^{\ell+1}_{U_2},\delta^{\ell+1}_{U^n},\delta^{\ell+1}_{H})|_{t=0}=(0,0,0,0),
\end{array}\right.
\end{eqnarray}
where
\beno
&&\delta\cF_1^{(\ell)}=\cF_1^{(\ell)}-\cF_1^{(\ell-1)}+\epsilon_0\f 1 {E^{(\ell)}\sqrt{E^{(\ell)}}}\Delta(\sqrt{E^{(\ell)}}U_1^{(\ell)})
-\epsilon_0\f 1 {E^{(\ell-1)}\sqrt{E^{(\ell-1)}}}\Delta(\sqrt{E^{(\ell-1)}}U_1^{(\ell)}),\\
&&\delta\cF_2^{(\ell)}=\cF_2^{(\ell)}-\cF_2^{(\ell-1)}+\epsilon_0\f 1 {E^{(\ell)}\sqrt{E^{(\ell)}}}\Delta(\sqrt{E^{(\ell)}}U_2^{(\ell)})
-\epsilon_0\f 1 {E^{(\ell-1)}\sqrt{E^{(\ell-1)}}}\Delta(\sqrt{E^{(\ell-1)}}U_2^{(\ell)}),\\
&&\delta\cF_3^{(\ell)}=\cF_3^{(\ell)}-\cF_3^{(\ell-1)}-\frac{1}{2E^{(\ell)}}\Delta H^{(\ell)}+\frac{1}{2E^{(\ell-1)}}\Delta H^{(\ell)}
-(2\delta_{u_i}^\ell-\delta_{w_i}^\ell)\p_iU^{n,(\ell)},\\
&&\delta\cF_4^{(\ell)}=F_7^{(\ell)}-F_7^{(\ell-1)}+\frac{1}{2E^{(\ell)}}\Delta U^{n,(\ell)}-\frac{1}{2E^{(\ell-1)}}\Delta U^{n,(\ell)}
+\delta_{w_i}^\ell\p_iH^{(\ell)}.
\eeno
From Proposition \ref{prop:nonlinear estimate-difference}, it is easy to see that
\beno
&&\|\delta\cF_i^{(\ell)}\|_{H^{s-4}}\le C\big(\|\delta_\RR^\ell\|_{H^{s-1}}+\|\delta_E^\ell\|_{H^{s-1}}
+\|(\delta^{\ell}_{U_1},\delta^{\ell}_{U_2},\delta^{\ell}_{U^n},\delta^{\ell}_{H})\|_{H^{s-3}}\big),\quad i=1,2,\\
&&\|\delta\cF_i^{(\ell)}\|_{H^{s-3}}\le C\big(\|\delta_\RR^\ell\|_{H^{s-1}}+\|\delta_E^\ell\|_{H^{s-1}}
+\|(\delta^{\ell}_{U_1},\delta^{\ell}_{U_2})\|_{H^{s-2}}+\|(\delta^{\ell}_{U^n},\delta^{\ell}_{H})\|_{H^{s-3}}\big),\quad i=3,4.
\eeno
Revisiting the proof of Theorem \ref{thm:linear}, we can obtain
\ben\label{eq:D_1}
D^{\ell+1}_1(t)&\le& C_\veps\int_0^tD^{\ell+1}_1(\tau)d\tau+C\sum_{i=1}^2\int_0^t\|\delta\cF_i^{(\ell}(\tau)\|_{H^{s-4}}^2d\tau\nonumber\\
&&+C\veps\sum_{i=3}^4\int_0^t\|\delta\cF_i^{(\ell}(\tau)\|_{H^{s-4}}^2d\tau\nonumber\\
&\le&C_\veps\int_0^tD^{\ell+1}_1(\tau)d\tau+C(t+\veps)\sup_{\tau \in [0,t]}D^{\ell}(\tau),
\een
where $D^\ell(t)$ is defined by
\beno
D^{\ell}(t)=D_1^{\ell}(t)+\|\delta_\RR^\ell(t)\|_{H^{s-1}}+\|\delta_E^\ell\|_{H^{s-1}},
\eeno
with
$
D^{\ell}_1(t)=\|(\delta^{\ell}_{U_1},\delta^{\ell}_{U_2},\delta^{\ell}_{U^n},\delta^{\ell}_{H})(t)\|_{H^{s-3}}
+\int_0^t\|(\delta^{\ell}_{U_1},\delta^{\ell}_{U_2})(\tau)\|_{H^{s-2}}^2d\tau.
$

On the other hand, we revisit the proof of Step 1 to find that
\ben\label{eq:D_2}
\|\delta_\RR^{\ell+1}\|_{H^{s-1}}+\|\delta_E^{\ell+1}\|_{H^{s-1}}\le C_0D^{\ell}_1(t)+Ct(\|\delta_\RR^{\ell}\|_{H^{s-1}}+\|\delta_E^{\ell}\|_{H^{s-1}}).
\een
For some small $\delta>0$ depending only on $C_0$, with $\veps$ and $T$ taken to be small enough, it follows from (\ref{eq:D_1}) and (\ref{eq:D_2}) that
\beno
\sup_{t\in [0,T]}(D_1^{\ell+1}(t)+\delta D_2^{\ell+1}(t))\le \f 12\sup_{t\in [0,T]}(D_1^{\ell}(t)+\delta D_2^{\ell}(t)),
\eeno
with $D_2^{\ell+1}=\|\delta_\RR^{\ell+1}\|_{H^{s-1}}+\|\delta_E^{\ell+1}\|_{H^{s-1}}$.
This implies that
$$(\RR^{(\ell)},U_1^{(\ell)}, U_2^{(\ell)}, U^{n, (\ell)},E^{(\ell)},H^{(\ell)},\Pi^{(\ell)},\widetilde{\RR}^{(\ell)},\widehat{\RR}^{(\ell)})$$
is a Cauchy sequence. More precisely, there exists the limit $(\RR,U_1, U_2, U^{n},E,H,\Pi,\widetilde{\RR},\widehat{\RR})$  such that
\beno
&&\RR^{(\ell)}\rightarrow\RR,\quad E^{(\ell)}\rightarrow E \quad \textrm{in}\quad L^\infty(0,T;H^{s-1});\\
&&U_1^{(\ell)}\rightarrow U_1,\quad  U_2^{(\ell)}\rightarrow U_2\quad \textrm{in}\quad L^\infty(0,T;H^{s-3})\cap L^2(0,T;H^{s-2});\\
&&U^{n, (\ell)}\rightarrow U^n,\quad H^{(\ell)}\rightarrow H,\quad \Pi^{(\ell)}\rightarrow\Pi \quad \textrm{in}\quad L^\infty(0,T;H^{s-3});\\
&&\widetilde{\RR}^{(\ell)}\rightarrow\widetilde{\RR}\quad \textrm{in}\quad L^\infty(0,T;H^{s-3}),\quad
\widehat{\RR}^{(\ell)}\rightarrow\widehat{\RR}\quad \textrm{in}\quad L^\infty(0,T;H^{s-2}).
\eeno
With the above information, it is easy to prove that $(\RR,U_1, U_2, U^{n},E,H,\Pi,\widetilde{\RR},\widehat{\RR})$
satisfies the system (\ref{eq:approximate system-velocity})-(\ref{eq:Pressure-app})
without the index $\ell$. In particular, we have
\begin{eqnarray}
&&\frac{\partial{\TR}}{\partial{t}}=\frac{\un}{E}\RRa\times\RRb+\frac{W_1}{\se}\RRa+\frac{W_2}{\se}\RRb,
\label{eq:limit-TR}\\
&&\Delta{\widehat\RR}-\HR=2H\TRa\times\TRb-\TR,\label{eq:limit-HR}\\
&&\Delta{\RR}-\RR=2H\HRa\times\HRb-\TR,\label{eq:limit-R}\\
&&\Delta{E}-2E=2(\RR_{\ab}\cdot\RR_{\ab}-\RR_{\alpha\alpha}\cdot\RR_{\beta\beta})-(\RRa\cdot\RRa+\RRb\cdot\RRb),
\label{eq:limit-E}\\
&&\frac{\partial{H}}{\partial{t}}=\frac{1}{2E}\Delta\un+\frac{W_1}{\se}
H_{\alpha}+\frac{W_2}{\se} H_{\beta}+\frac{\un}{2E^2}(L^2+2M^2+N^2)\label{eq:limit-H}\\
&&L=\RR_{\al\al}\cdot\nn,\quad M=\RR_{\al\be}\cdot\nn,\quad
N=\RR_{\be\be}\cdot\nn,
\quad\nn=\frac{\RR_{\alpha}\times\RR_{\beta}}{|\RR_{\alpha}\times\RR_{\beta}|}.
\end{eqnarray}
And, $(W_1,W_2)$ satisfies
\begin{eqnarray}\label{eq:limit-W}
\left\{
\begin{array}{l}
\Big(\frac{W_1}{\se}\Big)_{\alpha}-\Big(\frac{W_2}{\se}\Big)_{\beta}=\frac{\un(L-N)}{E},\\
\Big(\frac{W_1}{\se}\Big)_{\beta}+\Big(\frac{W_2}{\se}\Big)_{\alpha}=\frac{2\un M}{E}.
\end{array}\right.
\end{eqnarray}

It remains to show that the solution of the limit system is a solution of the original system.
For this purpose, it suffices to prove the following relations:
\ben\label{eq:relation}
\RR=\TR=\HR, ~~E=\RRa\cdot\RRa=\RRb\cdot\RRb,~~\RRa\cdot\RRb=0,\quad H=\frac{L+N}{2E}.
\een
And the incompressible condition follows easily from (\ref{eq:relation}).
As the proof is very complicated, it will be given in the following subsection.

\subsection{Consistency with the original system}

This subsection is devoted to proving (\ref{eq:relation}). Let us introduce some notations:
\beno
&&a_{11}=\RRa\cdot\RRa,\quad a_{12}=\RRa\cdot\RRb,\quad a_{22}=\RRb\cdot\RRb,\\
&&\widetilde{a}_{11}=\widetilde{\RR}_\al\cdot\widetilde{\RR}_\al,\quad \widetilde{a}_{12}
=\widetilde{\RR}_\al\cdot\widetilde{\RR}_\be,
\quad \widetilde{a}_{22}=\widetilde{\RR}_\be\cdot\widetilde{\RR}_\be,\\\
&&\ww=\frac{\un}{E}\RRa\times\RRb+\frac{W_1}{\se}\RRa+\frac{W_2}{\se}\RRb,\\
&&\overrightarrow{\ww}=(w_1,w_2),\quad w_i=\frac{W_i}{\se}(i=1,2).
\eeno
We set
\beno
\delta_R=\RR-\HR,\quad \delta_{\widetilde{R}}=\TR-\HR,\quad \delta_E^1=E-a_{11},\quad \delta_E^2=E-a_{22},\quad
\delta_a=\widetilde{a}_{11}-\widetilde{a}_{22}.
\eeno
In what follows, we denote by $\FF$ some operator bounded in $H^k(\RT^2)(0\le k\le 1)$,
which may be different from line to line. For example,
\beno
\|\FF(u,v)\|_{H^k}\le C(\|u\|_{H^k}+\|v\|_{H^k}).
\eeno

We get by using (\ref{eq:limit-TR}) and (\ref{eq:limit-W}) that
\begin{eqnarray}
&&(\tilde{a}_{11}-\tilde{a}_{22})_t
=\ww_{\alpha}\cdot\TRa-\ww_{\beta}\cdot\TRb\nonumber\\
&&=\ww_{\alpha}\cdot\RRa-\ww_{\beta}\cdot\RRb+\FF\big(\na(\RR-\TR)\big)\nonumber\\
&&=\frac{\un(L-N)}{2E}(2|\RRa\times\RRb|-a_{11}-a_{22})
+\frac{1}{2}(w_{1\alpha}+w_{2\beta})(a_{11}-a_{22})\nonumber\\
&&\quad+w_1(a_{11}-a_{22})_{\alpha}+w_2(a_{11}-a_{22})_{\beta}
+(w_{1\alpha}-w_{1\beta})a_{12}
+\FF\big(\na\delta_R,\na\delta_{\widetilde{R}}\big).\nonumber
\end{eqnarray}
Similarly, we have
\begin{eqnarray}
(\tilde{a}_{12})_t
&=&\ww_{\alpha}\cdot\TRb+\ww_{\beta}\cdot\TRa\nonumber\\
&=&\ww_{\alpha}\cdot\RRb+\ww_{\beta}\cdot\RRa+\FF\big(\na(\RR-\TR)\big)\nonumber\\
&=&\frac{M\un}{E}(a_{11}+a_{22}-2|\RRa\times\RRb|)
+\frac{1}{2}(w_{1\beta}-w_{2\alpha})(a_{11}-a_{22})\nonumber\\
&&+w_1a_{12\alpha}+w_2a_{12\beta}+(w_{1\alpha}+w_{2\beta})a_{12}+\FF\big(\na\delta_R,\na\delta_{\widetilde{R}}\big).\nonumber
\end{eqnarray}
Noting
\begin{eqnarray}
a_{11}+a_{22}-2|\RRa\times\RRb|
&=&\frac{(a_{11}+a_{22})^2-4(a_{11}a_{22}-a_{12}^2)}{a_{11}+a_{22}+2|\RRa\times\RRb|}\nonumber\\
&=&\frac{(a_{11}-a_{22})^2-4a_{12}^2}{a_{11}+a_{22}+2|\RRa\times\RRb|},\nonumber
\end{eqnarray}
we find that
\ben
&&\p_t\delta_{a}+\overrightarrow{\ww}\cdot\na\delta_{a}=\FF\big(\na\delta_R,\na\delta_{\widetilde{R}},\na^2\delta_R,\na^2\delta_{\widetilde{R}},\delta_a,\widetilde{a}_{12}\big),\label{eq:con-ta11}\\
&&\p_t\tilde{a}_{12}+\overrightarrow{\ww}\cdot\na \widetilde{a}_{12}=\FF\big(\na\delta_R,\na\delta_{\widetilde{R}},\na^2\delta_R,\na^2\delta_{\widetilde{R}},\delta_a,\widetilde{a}_{12}\big).\label{eq:con-ta12}
\een
By (\ref{eq:limit-HR}) and (\ref{eq:limit-R}), we have
\ben\label{eq:con-HR-R}
(\Delta-1)(\HR-\RR)=2H(\TRa\times\TRb-\HRa\times\HRb)=\FF\big(\na\delta_{\widetilde{R}}\big).
\een
And by (\ref{eq:limit-HR}) and (\ref{eq:limit-E}),
\beno
(\Delta-1)(E-a_{11})=-2\big(\RRa\cdot(2H\HRa\times\HRb)\big)_{\alpha}+E-a_{22}-2(\RRa\cdot(\RR-\TR))_{\alpha},
\eeno
which implies that
\ben\label{eq:con-E-a11}
(\Delta-2)\delta_E^1=\FF\big(\na\delta_R,\na\delta_{\widetilde{R}},\delta_a\big).
\een
Similarly, we have
\ben\label{eq:con-E-a22}
(\Delta-2)(E-a_{22})&=&-2\big(\RRb\cdot(2H\HRa\times\HRb)\big)_{\beta}+a_{22}-a_{11}-2(\RRb\cdot(\RR-\TR))_{\beta}\nonumber\\
&=&\FF\big(\na\delta_R,\na\delta_{\widetilde{R}},\delta_a\big),\\
(\Delta-1)a_{12}&=&2H\big((\HRa\times\HR_{\ab})\cdot\RRb-(\HRa\times\RR_{\ab})\cdot\HRb\big)\nonumber\\
&&+2H\big((\HR_{\ab}\times\HRb)\cdot\RRa-(\RR_{\ab}\times\HRb)\cdot\HRa\big)\nonumber\\
&&-a_{12}+\FF(\delta_R,\nabla\delta_R,\delta_{\widetilde{R}},\na\delta_{\widetilde{R}})\nonumber\\
&=&\sum_{k=0}^2\FF\big(\na^k\delta_R,\na^k\delta_{\widetilde{R}},\widetilde{a}_{12}\big).\label{eq:con-a12}
\een
The following facts will be used frequently:
\ben
&&|\RRa\times\RRb|-E=\sqrt{\det(a_{ij})}-E=\FF\big(\delta_E^1,\delta_E^2,a_{12}\big),\label{fact-1}\\
&&\Gamma_{11}^1,\Gamma_{12}^2,\Gamma_{21}^2,-\Gamma_{22}^1=\frac{\Ea}{2E}+\FF\big(\delta_E^1,\delta_E^2,\na\delta_E^1,\na\delta_E^2, a_{12}\big),\label{fact-2}\\
&&\Gamma_{12}^1,\Gamma_{21}^1,\Gamma_{22}^2,-\Gamma_{11}^2=\frac{\Eb}{2E}+\FF\big(\delta_E^1,\delta_E^2,\na\delta_E^1,\na\delta_E^2, a_{12}\big).\label{fact-3}\\
&&(\Delta\RR)_{\gamma}=(2H\RRa\times\RRb)_{\gamma}+\FF\big(\na\delta_R,\na\delta_{\widetilde{R}},\na^2\delta_R,\na^2\delta_{\widetilde{R}}\big),
\quad\gamma=\al,\be.\label{fact-4}
\een
Indeed, we have
\beno
\Gamma_{12}^2&=&\f 1 {|\RRa\times\RRb|}\big(-\f {a_{12}}2\p_\be a_{11}+\f {a_{11}} 2\p_\al a_{22}\big)\\
&=&\frac{\Ea}{2E}+\FF\big(\delta_E^1,\delta_E^2,\na\delta_E^2, a_{12}\big),
\eeno
and the others can be deduced similarly. For the last fact, we have by (\ref{eq:limit-R}) that
\ben
\Delta\RR&=&2H\HR_\al\times \HR_\be-\TR+\RR\nonumber\\
&=&2H\RR_\al\times \RR_\be+2H(\HR_\al\times \HR_\be-\RR_\al\times \RR_\be)-\TR+\RR,\label{eq:Delta-RR}
\een
thus (\ref{fact-4}) follows easily.

To proceed, we also need the following lemma.
\begin{Lemma}\label{lem:second fundamental}
For $\gamma=\al,\be$, it holds that
\ben
&&(2EH-L-N)_\gamma=\sum_{k=0}^2\sum_{j=0}^1\FF\big(\na^k\delta_R,\na^k\delta_{\widetilde{R}},\na^{j}\delta_E^1,\na^j\delta_E^2,\na^ja_{12}\big),\label{eq:sf-1}\\
&&LN-M^2=\frac{1}{2}\Big(-\Delta{E}+\frac{\Ea^2+\Eb^2}{E}\Big)+\sum_{k=0}^2\FF\big(\na^k\delta_E^1,\na^k\delta_E^2,\na^k a_{12}\big),\label{eq:sf-2}\\
&&N_{\alpha}=M_{\beta}+H\Ea+\FF\big(\delta_E^1,\delta_E^2,\na\delta_E^1,\na\delta_E^2, a_{12}\big),\label{eq:sf-3}\\
&&L_{\beta}=M_{\alpha}+H\Eb+\FF\big(\delta_E^1,\delta_E^2,\na\delta_E^1,\na\delta_E^2, a_{12}\big).\label{eq:sf-4}
\een
\end{Lemma}

\no{\bf Proof.}\,First of all, a direct calculation gives
\beno
2EH-(L+N)&=&\frac{E}{|\HRa\times\HRb|^2}\Delta\RR\cdot(\HRa\times\HRb)-\frac{1}{|\RRa\times\RRb|}\Delta\RR\cdot(\RRa\times\RRb)\\
&=&\Delta\RR\cdot\Big(\frac{E\HRa\times\HRb}{|\HRa\times\HRb|^2}-\frac{\RRa\times\RRb}{|\RRa\times\RRb|}\Big),
\eeno
which implies (\ref{eq:sf-1}). From the Gauss equation, we infer that
\beno
{LN-M^2}&=&\frac{1}{a_{11}a_{22}-a_{12}^2}
\left(\Bigg|
\begin{array}{ccc}
a_{12\ab}-\frac{1}{2}a_{11\beta\beta}-\frac{1}{2}a_{22\alpha\alpha}
& \frac{1}{2}a_{11\alpha} & a_{12\alpha}-\frac{1}{2}a_{11\beta}\\
a_{12\beta}-\frac{1}{2}a_{22\alpha} & a_{11} & a_{12}\\
\frac{1}{2}a_{22\beta}&a_{12}&a_{22}
\end{array}\Bigg|\right.\\
&&-\left.\Bigg|
\begin{array}{ccc}
0&\frac{a_{11\beta}}{2}&\frac{a_{22\alpha}}{2}\\
\frac{a_{11\beta}}{2}&a_{11}&a_{12}\\
\frac{a_{22\alpha}}{2}&a_{12}&a_{22}
\end{array}
\Bigg|\right),
\eeno
which implies (\ref{eq:sf-2}). The Codazzi equation $b_{\alpha\beta,\gamma}=b_{\alpha\gamma,\beta}$ implies that
\beno
&&\frac{\p b_{11}}{\partial x^2}-\Gamma_{12}^1b_{11}-\Gamma_{12}^2b_{12}-\Gamma_{12}^1b_{11}-\Gamma_{12}^2b_{21}\\
&&\quad=\frac{\p b_{12}}{\partial x^1}-\Gamma_{21}^1b_{11}-\Gamma_{21}^2b_{12}-\Gamma_{11}^1b_{12}-\Gamma_{11}^2b_{22},\\
&&\frac{\p b_{21}}{\partial x^2}-\Gamma_{12}^1b_{21}-\Gamma_{12}^2b_{22}-\Gamma_{22}^1b_{11}-\Gamma_{22}^2b_{21}\\
&&\quad=\frac{\p b_{22}}{\partial x^1}-\Gamma_{21}^1b_{21}-\Gamma_{21}^2b_{22}-\Gamma_{21}^1b_{12}-\Gamma_{21}^2b_{22},
\eeno
where $b_{11}=L$, $b_{12}=b_{21}=M$, and $b_{22}=N$. Then (\ref{eq:sf-3})-(\ref{eq:sf-4}) follow easily
from (\ref{fact-2}) and (\ref{fact-3}).
The proof is completed.\endproof

In the following, we calculate $\TR-\HR$. By (\ref{eq:limit-TR}) and (\ref{eq:limit-HR}), we have
\ben
(\Delta-1)(\TR-\HR)_t=\Delta \ww-(2H\TRa\times\TRb)_t.\label{eq:TR-HR-1}
\een
Direct calculations yield that
\begin{eqnarray}
\Delta \ww&=&\Delta(\frac{\un}{E}\RRa\times\RRb+w_1\RRa+w_2\RRb)\nonumber\\
&=&\Delta\Big(\frac{\un}{E}\Big)\RRa\times\RRb+\frac{\un}{E}
(\Delta\RR)_{\alpha}\times\RRb+\frac{\un}{E}\RRa\times(\Delta\RR)_{\beta}\nonumber\\
&&+\Delta{w_1}\RRa+\Delta{w_2}\RRb+w_1(\Delta\RR)_{\alpha}+w_2(\Delta\RR)_{\beta}
+2\frac{\un}{E}(\RR_{\alpha\alpha}-\RR_{\beta\beta})\times\RR_{\alpha\beta}\nonumber\\
&&+2\Big(\frac{\un}{E}\Big)_{\alpha}\RR_{\alpha\alpha}\times\RRb+2\Big(\frac{\un}{E}\Big)_{\alpha}\RRa\times\RR_{\ab}\nonumber\\
&&+2\Big(\frac{\un}{E}\Big)_{\beta}\RR_{\ab}\times\RRb+2\Big(\frac{\un}{E}\Big)_{\beta}\RRa\times\RR_{\beta\beta}\nonumber\\
&&+2w_{1\alpha}\RR_{\alpha\alpha}+2w_{2\beta}\RR_{\beta\beta}+2(w_{1\beta}+w_{2\alpha})\RR_{\ab},
\end{eqnarray}
and by (\ref{eq:limit-TR}),
\begin{eqnarray}
&&\Big(2H(\TRa\times\TRb)\Big)_t\nonumber\\
&&=\Big(\frac{\Delta\un}{E}+2w_1H_{\alpha}
+2w_2H_{\beta}+\frac{\un}{E^2}(L^2+2M^2+N^2)\Big)\RRa\times\RRb\nonumber\\
&&\quad+2H\Big(\frac{\un}{E}\RRa\times\RRb+w_1\RRa+w_2\RRb\Big)_{\alpha}\times\TRb\nonumber\\
&&\quad+2H\TRa\times\Big(\frac{\un}{E}\RRa\times\RRb+w_1\RRa+w_2\RRb\Big)_{\beta}\nonumber\\
&&=\Big(\frac{\Delta\un}{E}+2w_1H_{\alpha}
+2w_2H_{\beta}+\frac{\un}{E^2}(L^2+2M^2+N^2)\Big)\RRa\times\RRb\nonumber\\
&&\quad+2H\Big(\frac{\un}{E}\RRa\times\RRb+w_1\RRa+w_2\RRb\Big)_{\alpha}\times\RRb\nonumber\\
&&\quad+2H\RRa\times\Big(\frac{\un}{E}\RRa\times\RRb+w_1\RRa+w_2\RRb\Big)_{\beta}+\FF\big(\na\delta_R,\na\delta_{\widetilde{R}}\big)\nonumber\\
&&=\Big(\frac{\Delta\un}{E}+\frac{\un}{E^2}(L^2+2M^2+N^2)\Big)\RRa\times\RRb
+w_1(2H\RRa\times\RRb)_{\alpha}\nonumber\\
&&\quad+w_2(2H\RRa\times\RRb)_{\beta}+2H(w_{1\alpha}+w_{2\beta})\RRa\times\RRb\nonumber\\
&&\quad+2H\Big(\frac{\un}{E}\Big)_{\alpha}(\RRa\times\RRb)\times\RRb+2H\frac{\un}{E}\big(\RRa\times\RRb\big)_{\alpha}\times\RRb\nonumber\\
&&\quad+2H\Big(\frac{\un}{E}\Big)_{\beta}\RRa\times(\RRa\times\RRb)+2H\frac{\un}{E}\RRa\times\big(\RRa\times\RRb\big)_{\beta}
+\FF\big(\na\delta_R,\na\delta_{\widetilde{R}}\big).
\end{eqnarray}
From (\ref{fact-2})-(\ref{fact-3}), it follows that
\begin{eqnarray}
\RR_{\alpha\alpha}&=&\Gamma_{11}^1\RRa+\Gamma_{11}^2\RRb+L\nn\nonumber\\
&=&\frac{\Ea}{2E}\RRa-\frac{\Eb}{2E}\RRb+L\nn+\FF\big(\delta_E^1,\delta_E^2,\na\delta_E^1,\na\delta_E^2, a_{12}\big),\label{eq:TR-HR-4}\\
\RR_{\alpha\beta}&=&\Gamma_{12}^1\RRa+\Gamma_{12}^2\RRb+M\nn\nonumber\\
&=&\frac{\Eb}{2E}\RRa+\frac{\Ea}{2E}\RRb+M\nn+\FF\big(\delta_E^1,\delta_E^2,\na\delta_E^1,\na\delta_E^2, a_{12}\big),\label{eq:TR-HR-5}\\
\RR_{\beta\beta}&=&\Gamma_{22}^1\RRa+\Gamma_{22}^2\RRb+N\nn\nonumber\\
&=&-\frac{\Ea}{2E}\RRa+\frac{\Eb}{2E}\RRb+N\nn
+\FF\big(\delta_E^1,\delta_E^2,\na\delta_E^1,\na\delta_E^2, a_{12}\big).\label{eq:TR-HR-6}
\end{eqnarray}
And by (\ref{fact-4}), it follows that
\begin{eqnarray}
&&(\Delta\RR)_{\alpha}\times\RRb\nonumber\\
&&=2H_{\alpha}(\RRa\times\RRb)\times\RRb+2H(\RRa\times\RRb)_{\alpha}\times\RRb
+\FF\big(\na\delta_R,\na\delta_{\widetilde{R}},\na^2\delta_R,\na^2\delta_{\widetilde{R}}\big)\nonumber\\
&&=-2EH_{\alpha}\RRa+2H(\RRa\times\RRb)_{\alpha}\times\RRb
+\FF\big(\na\delta_R,\na\delta_{\widetilde{R}},\na^2\delta_R,\na^2\delta_{\widetilde{R}},\delta_E,a_{12}\big),\\
&&\RRa\times(\Delta\RR)_{\beta}\nonumber\\
&&=2H_{\beta}\RRa\times(\RRa\times\RRb)+2H\RRa\times(\RRa\times\RRb)_{\beta}
+\FF\big(\na\delta_R,\na\delta_{\widetilde{R}},\na^2\delta_R,\na^2\delta_{\widetilde{R}}\big)\nonumber\\
&&=-2EH_{\beta}\RRb+2H\RRa\times(\RRa\times\RRb)_{\beta}
+\FF\big(\na\delta_R,\na\delta_{\widetilde{R}},\na^2\delta_R,\na^2\delta_{\widetilde{R}},\delta_E,a_{12}\big).\label{eq:TR-HR-8}
\end{eqnarray}

Summing up (\ref{eq:TR-HR-1})-(\ref{eq:TR-HR-8}), we obtain
\beno
&&(\Delta-1)(\TR-\HR)_t=\Delta \ww-(2H\TRa\times\TRb)_t\\
&&=\Delta\Big(\frac{\un}{E}\Big)E\nn-2\un(H_{\alpha}\RRa+H_{\beta}\RRb)+\Delta{w_1}\RRa+\Delta{w_2}\RRb\\
&&\quad+2\Big(\frac{\un}{E}\Big)_{\alpha}(\Ea\nn-L\RRa-M\RRb)+2\Big(\frac{\un}{E}\Big)_{\beta}(\Eb\nn-M\RRa-N\RRb)\\
&&\quad+2(w_{1\alpha}\RR_{\alpha\alpha}+w_{2\beta}\RR_{\beta\beta})-(w_{1\alpha}+w_{2\beta})(2H\RRa\times\RRb)\\
&&\quad+2\frac{\un}{E}(\RR_{\alpha\alpha}-\RR_{\beta\beta})\times\RR_{\alpha\beta}+\frac{4M\un}{E}\RR_{\alpha\beta}
-\Big({\Delta\un}+\frac{\un}{E}(L^2+2M^2+N^2)\Big)\nn\\
&&\quad+2EH\Big(\frac{\un}{E}\Big)_{\alpha}\RRa
+2EH\Big(\frac{\un}{E}\Big)_{\beta}\RRb\\
&&\quad+\FF\big(\na\delta_R,\na\delta_{\widetilde{R}},\na^2\delta_R,\na^2\delta_{\widetilde{R}},\delta_E^1,\delta_E^2,\na\delta_E^1,\na\delta_E^2, a_{12}\big).
\eeno
On the other hand, we can get by (\ref{eq:TR-HR-4})-(\ref{eq:TR-HR-6}) that
\beno
&&(\RR_{\alpha\alpha}-\RR_{\beta\beta})\times\RR_{\alpha\beta}\\
&&=\frac{\Ea^2+\Eb^2}{2E}\nn-\frac{M\Ea}{E}\RRb-\frac{M\Eb}{E}\RRa+\frac{(L-N)\Eb}{2E}\RRb-\frac{(L-N)\Ea}{2E}\RRa\\
&&\quad+\FF\big(\delta_E^1,\delta_E^2,\na\delta_E^1,\na\delta_E^2, a_{12}\big),
\eeno
and by (\ref{eq:Delta-RR}) and (\ref{eq:limit-W}),
\beno
&&2(w_{1\alpha}\RR_{\alpha\alpha}+w_{2\beta}\RR_{\beta\beta})-(w_{1\alpha}+w_{2\beta})(2H\RRa\times\RRb)\\
&&=(w_{1\alpha}-w_{2\beta})(\RR_{\alpha\alpha}-\RR_{\beta\beta})+\FF\big(\delta_R,\delta_{\widetilde{R}}, \na\delta_R,\na\delta_{\widetilde{R}}\big)\\
&&=\frac{\un(L-N)}{E}(\RR_{\alpha\alpha}-\RR_{\beta\beta})+\FF\big(\delta_R,\delta_{\widetilde{R}},\na\delta_R,\na\delta_{\widetilde{R}}\big)\\
&&=\frac{\un(L-N)}{E}\Big(\frac{\Ea}{E}\RRa-\frac{\Eb}{E}\RRb+(L-N)\nn\Big)\\
&&\quad+\FF\big(\delta_R,\delta_{\widetilde{R}},\na\delta_R,\na\delta_{\widetilde{R}},\delta_E^1,\delta_E^2,\na\delta_E^1,\na\delta_E^2, a_{12}\big).
\eeno
Thus, we arrive at
\beno
&&(\Delta-1)(\TR-\HR)_t\\
&&=\Delta\Big(\frac{\un}{E}\Big)E\nn-2\un(H_{\alpha}\RRa+H_{\beta}\RRb)+\Delta{w_1}\RRa+\Delta{w_2}\RRb\\
&&\quad+2\Big(\frac{\un}{E}\Big)_{\alpha}(\Ea\nn-L\RRa-M\RRb)+2\Big(\frac{\un}{E}\Big)_{\beta}(\Eb\nn-M\RRa-N\RRb)\\
&&\quad+\frac{\un}{E}\Big(\frac{\Ea^2+\Eb^2}{E}\nn-\frac{2M\Ea}{E}\RRb-\frac{2M\Eb}{E}\RRa+\frac{(L-N)\Eb}{E}\RRb-\frac{(L-N)\Ea}{E}\RRa\Big)\\
&&\quad+\frac{\un(L-N)}{E}\Big(\frac{\Ea}{E}\RRa-\frac{\Eb}{E}\RRb+(L-N)\nn\Big)
+\frac{4M\un}{E}\Big(\frac{\Eb}{2E}\RRa+\frac{\Ea}{2E}\RRb+M\nn\Big)\\
&&\quad-\Big({\Delta\un}+\frac{\un}{E}(L^2+2M^2+N^2)\Big)\nn
+2EH\Big(\frac{\un}{E}\Big)_{\alpha}\RRa
+2EH\Big(\frac{\un}{E}\Big)_{\beta}\RRb\\
&&\quad+\FF\big(\delta_R,\delta_{\widetilde{R}},\na\delta_R,\na\delta_{\widetilde{R}},\na^2\delta_R,\na^2\delta_{\widetilde{R}},\delta_E^1,\delta_E^2,\na\delta_E^1,\na\delta_E^2, a_{12}\big)\\
&&=\frac{\un}{E}\Big(2(M^2-LN)-\Delta{E}+\frac{\Ea^2+\Eb^2}{E}\Big)\nn\\
&&\quad+\Big(-2\un{H_{\alpha}}+\Delta{w_1}-(L-N)\big(\frac{\un}{E}\big)_{\alpha}-2M\big(\frac{\un}{E}\big)_{\beta}\Big)\RRa\\
&&\quad+\Big(-2\un{H_{\beta}}+\Delta{w_2}-2M\big(\frac{\un}{E}\big)_{\alpha}+(L-N)\big(\frac{\un}{E}\big)_{\beta}\Big)\RRb\\
&&\quad+\FF\big(\delta_R,\delta_{\widetilde{R}},\na\delta_R,\na\delta_{\widetilde{R}},\na^2\delta_R,\na^2\delta_{\widetilde{R}},\delta_E^1,\delta_E^2,\na\delta_E^1,\na\delta_E^2, a_{12}\big).
\eeno
And thanks  to (\ref{eq:limit-W}), we have
\beno
\Delta{w_1}=\Big(\frac{\un(L-N)}{E}\Big)_{\alpha}+\Big(\frac{2M\un}{E}\Big)_{\beta},\\
\Delta{w_2}=\Big(\frac{2M\un}{E}\Big)_{\alpha}-\Big(\frac{\un(L-N)}{E}\Big)_{\beta},
\eeno
which together with Lemma \ref{lem:second fundamental} implies that
\begin{eqnarray}\label{eq:con-TR-HR}
&&(\Delta-1)(\TR-\HR)_t\nonumber\\
&&=\frac{\un}{E}\Big(2(M^2-LN)-\Delta{E}+\frac{\Ea^2+\Eb^2}{E}\Big)\nn\nonumber\\
&&\quad+\frac{\un}{E}\big(L+N-2EH\big)_{\alpha}\RRa
+\frac{\un}{E}\big(L+N-2EH\big)_{\beta}\RRb\nonumber\\
&&\quad+\sum_{k=0}^2\FF\big(\na^k\delta_R,\na^k\delta_{\widetilde{R}},\na^k\delta_E^1,\na^k\delta_E^2, \na^k a_{12}\big)\nonumber\\
&&=\sum_{k=0}^2\FF\big(\na^k\delta_R,\na^k\delta_{\widetilde{R}},\na^k\delta_E^1,\na^k\delta_E^2, \na^k a_{12}\big).
\end{eqnarray}

Now we are position to prove (\ref{eq:relation}).
Firstly, by (\ref{eq:R-app-tilde}) and the fact that the initial surface is
parameterized by the isothermal coordinates and
(\ref{eq:R-elliptic equation}) holds for $t=0$, we know that all the relations in (\ref{eq:relation}) hold for $t=0$.
Hence,
\beno
\delta_a(0)=0,\quad \widetilde{a}_{12}(0)=0,\quad \delta_{\widetilde{R}}(0)=0.
\eeno
Taking the $L^2$ energy estimate to (\ref{eq:con-ta11}) and (\ref{eq:con-ta12}), we obtain
\beno
&&\|\delta_a(t)\|_{L^2}\le C\int_0^t\big(\|\delta_a\|_{L^2}+\|\widetilde{a}_{12}\|_{L^2}+\|\delta_R\|_{H^2}+\|\delta_{\widetilde{R}}\|_{H^2}\big)d\tau,\\
&&\|\widetilde{a}_{12}(t)\|_{L^2}\le C\int_0^t\big(\|\delta_a\|_{L^2}+\|\widetilde{a}_{12}\|_{L^2}+\|\delta_R\|_{H^2}+\|\delta_{\widetilde{R}}\|_{H^2}\big)d\tau.
\eeno
Using the elliptic estimate, we deduce from (\ref{eq:con-HR-R})-(\ref{eq:con-a12}) that
\beno
&&\|\delta_R\|_{H^2}\le C\|\delta_{\widetilde{R}}\|_{H^1},\\
&&\|(\delta_E^1, \delta_E^2,a_{12})\|_{H^2}\le C\big(\|\widetilde{a}_{12}\|_{L^2}+\|\delta_R\|_{H^1}+\|\delta_{\widetilde{R}}\|_{H^1}\big),
\eeno
and from (\ref{eq:con-TR-HR}), it follows that
\beno
\|\delta_{\widetilde{R}}(t)\|_{H^2}\le C\int_0^t\big(\|\delta_a\|_{L^2}+\|\widetilde{a}_{12}\|_{L^2}+\|\delta_R\|_{H^2}+\|\delta_{\widetilde{R}}\|_{H^2}
+\|(\delta_E^1, \delta_E^2,a_{12})\|_{H^2}\big)d\tau.
\eeno
Thus, we obtain
\beno
&&\|\delta_a(t)\|_{L^2}+\|\widetilde{a}_{12}(t)\|_{L^2}+\|\delta_{\widetilde{R}}(t)\|_{H^2}\\
&&\le C\int_0^t\big(\|\delta_a(\tau)\|_{L^2}+\|\widetilde{a}_{12}(\tau)\|_{L^2}+\|\delta_{\widetilde{R}}(\tau)\|_{H^2}\big)d\tau,
\eeno
which implies (\ref{eq:relation}) by Gronwall's inequality.

\subsection{Remark on the general case}

In this subsection, we describe how to adapt our method to deal with
the case in which the surface is parameterized by a finite number of
isothermal coordinates. Assume that we need $N$ local chart to
parameterize the initial surface $S_0=\cup_{i=1}^NS_0^i$ where each
$S_0^i$ is open and parameterized by isothermal coordinates:
\beno \RR_0^i(x_1,x_2):\Om^i\longrightarrow S_0^i,\quad 1\le i\le N.
\eeno
Let $\{\psi^i\}_{1\le i\le N}$ be a partition of the unit
subordinate to $\{S_0^i\}_{1\le i\le N}$; that is, \beno
\sum_{i=1}^N\psi^i=1,\quad \textrm{supp}\psi^i\subset S_0^i. \eeno
At each local chart, $\RR^i$ is defined by \beno
&&\frac{\partial\RR^i}{\partial t}=v^{n}\nn_i+W_1^i\ta_i+W_2^i\tb_i,
\eeno where $(W^i_1,W^i_2)$ is defined by \beno \left\{
\begin{array}{l}
\Big(\frac{W_1^i}{\sqrt{E_i}}\Big)_{\alpha}-\Big(\frac{W_2^i}{\sqrt{E_i}}\Big)_{\beta}=\frac{\un(L_i-N_i)}{E_i},\\
\Big(\frac{W_1^i}{\sqrt{E_i}}\Big)_{\beta}+\Big(\frac{W_2^i}{\sqrt{E_i}}\Big)_{\alpha}=\frac{2\un M_i}{E_i}.
\end{array}\right.
\eeno
While, $(v,H,\Pi)$ is determined by solving the following system:
\beno
&&\frac{\partial{\vv}}{\partial{t}}=(-\Pi{a^{\alpha\beta}\aaa_{\alpha}})_{,\beta}+2\varepsilon_0
(S^{\al\be}\aaa_{\alpha})_{,\beta}
-\frac{1}{2}\Big(\Delta_{\Gamma}{H}+H(b^{\alpha}_{\beta}b^{\beta}_{\alpha}-2H^2)\Big)\nn,\\
&&v^{\alpha}_{,\alpha}=2Hv^{n},\\
&&2\frac{\partial{H}}{\partial{t}}=a^{\al\be}v^n_{,\al\be}+v^nb^{\alpha}_{\beta}b^{\beta}_{\alpha}+2v^{\alpha}H_{,\alpha},
\eeno see Section 2 for some notations. As the above equations
are coordinate-invariant, $(v,H,\Pi)$ does not depend on the choice
of coordinates. In this case, the energy functional is given by
\beno
\cE(t)=\|(-\Delta_{S_t})^kv^T\|_{L^2(S_t)}^2+\|(-\Delta_{S_t})^{k}v^n\|_{L^2(S_t)}^2+\|(-\Delta_{S_t})^{k}H\|_{L^2(S_t)}^2,
\eeno
where $v^T$ is the tangential component of the velocity, and
$\Delta_{S_t}$ is the Laplace-Beltrami operator on the surface $S_t$
at time $t$. In the isothermal coordinates, $\Delta_{S_t}=\f 1 E\Delta$.
Then, as in section 4, we can obtain a uniform estimate for
$\cE(t)$. Let $\{\phi^i(t,x_1,x_2)\}_{1\le i\le N}$ be a partition
of the unit on $S_t$ given by \beno \phi^i(t,x_1,x_2)=\f
{\psi^i(R^i_0(x_1,x_2))}
{\sum_{j=1}^N\widetilde{\psi}^j(t,R^i(t,x_1,x_2))}, \quad
\widetilde{\psi}^i(t,X)=\psi^i(R^i_0\circ(R^i(t))^{-1}(X)). \eeno
Indeed, we have \beno \f d
{dt}\|(-\Delta_{S_t})^{k}v^n\|_{L^2(S_t)}^2&=&\sum_{i=1}^N\p_t\int\phi^i\big|(-\f
1 {E_i}\Delta)^kv^n\big|^2E_idx_1dx_2\\
&=&-\sum_{i=1}^N\int\phi^i(-\f 1 {E_i}\Delta)^kv^n(-\f 1 {E_i}\Delta)^{k+1}HE_idx_1dx_2+L.W.T.\\
&=&-\int_{S_t}(-\Delta_{S_t})^kv^n(-\Delta_{S_t})^{k+1}H dS_t+L.W.T..
\eeno
And, similarly,
\beno
\f d {dt}\|(-\Delta_{S_t})^{k}H\|_{L^2(S_t)}^2&=&\sum_{i=1}^N\int\phi^i(-\f 1 {E_i}\Delta)^kH(-\f 1 {E_i}\Delta)^{k+1}v^nE_idx_1dx_2+L.W.T.\\
&=&\int_{S_t}(-\Delta_{S_t})^kH(-\Delta_{S_t})^{k+1}v^ndS_t+L.W.T.,
\eeno
where $L.W.T.$ denotes the lower-order terms. Thus, we have
\beno
\f d {dt}\big(\|(-\Delta_{S_t})^{k}v^n\|_{L^2(S_t)}^2+\|(-\Delta_{S_t})^{k}H\|_{L^2(S_t)}^2\big)=L.W.T.
\eeno

\setcounter{equation}{0}
\section{Appendix}

\subsection{Derivations of the equation (\ref{eq:velocity-2}) and the energy law}
In this subsection, we give the derivations of the equation
(\ref{eq:velocity-2}) and the energy law in the case
$B_{\alpha\beta}=Ba_{\alpha\beta}$. The reader can also find a short
version of the derivation in \cite{Hu-Zhang}.

Firstly, we have \beno
M^{\alpha\beta}=C_1^{\alpha\beta\gamma\delta}(Ba_{\gamma\delta}-b_{\gamma\delta})
=2(k_1B-(k_1-\varepsilon_1)H)a^{\alpha\beta}+2\varepsilon_1b^{\alpha\beta}.
\eeno We get by (\ref{eq:Guass-Weiggarten-Codazzi}) that \beno
(M^{\alpha\mu}b_{\mu}^{\beta}\aaa_{\beta})_{,\alpha}+(q^{\alpha}\nn)_{\alpha}
&=&(M^{\alpha\mu}b_{\mu}^{\beta})_{,\alpha}\aab+M^{\alpha\mu}b_{\mu}^{\beta}
b_{\alpha\beta}\nn+M^{\alpha\beta}_{,\alpha\beta}\nn-M^{\alpha\gamma}_{\gamma}b^{\beta}_{\alpha}\aab\nonumber\\
&=&M^{\alpha\mu}b^{\beta}_{\mu,\alpha}\aab+(M^{\alpha\mu}b^{\beta}_{\mu}
b_{\alpha\beta}+M^{\alpha\beta}_{,\alpha\beta})\nn.
\eeno
Using $b_{\alpha\beta,\gamma}=b_{\alpha\gamma,\beta}$, we infer that
\begin{eqnarray}
&&a^{\alpha\mu}b^{\beta}_{\mu,\alpha}=b^{\alpha\beta}_{,\alpha}=a^{\beta\mu}b_{\mu,\alpha}^{\alpha}
=a^{\beta\mu}b^{\alpha}_{\alpha,\mu}=2a^{\alpha\beta}H_{,\alpha},\nonumber\\
&&\frac{1}{2}a^{\alpha\beta}(b^{\gamma\delta}b_{\gamma\delta})_{,\alpha}
=a^{\alpha\beta}b^{\gamma\delta}b_{\gamma\delta,\alpha}
=a^{\alpha\beta}b^{\gamma\delta}b_{\gamma\alpha,\delta}
=b^{\gamma\delta}b^{\beta}_{\gamma,\delta}
=b^{\alpha\mu}b_{\mu,\alpha}^{\beta},\nonumber
\end{eqnarray}
from which, we can deduce that
\begin{eqnarray}
M^{\alpha\mu}b^{\beta}_{\mu,\alpha}&=&2
(k_1B-(k_1-\varepsilon_1)H)b^{\alpha\beta}_{,\beta}+2\varepsilon_1b^{\alpha\mu}
b^{\beta}_{\mu,\al}\nonumber\\
&=&-4k_1Ha^{\alpha\beta}B_{,\alpha}-a^{\alpha\beta}\big(2(k_1-
\varepsilon_1)H^2+\varepsilon_1b^{\gamma\delta}b_{\gamma\delta}-4k_1HB\big)_{,\alpha},\nonumber\\
M^{\alpha\beta}_{,\alpha\beta}&=&2k_1a^{\alpha\beta}B_{,\alpha\beta}
-2(k_1+\varepsilon_1)a^{\alpha\beta}H_{,\alpha\beta}.\nonumber
\end{eqnarray}

Let $K$ be the Gaussian curvature. It is easy to see that
\begin{equation*}
b^{\alpha\beta}b_{\alpha\beta}=4H^2-2K,\quad b^{\alpha\beta}b^{\gamma}_{\beta}b_{\alpha\gamma}=2H(4H^2-3K),
\end{equation*}
which implies that
\beno
M^{\alpha\mu}b^{\beta}_{\mu}b_{\alpha\beta}=4k_1(B-H)(2H^2-K)-8\varepsilon_1H(H^2-K).
\eeno
Let $P\eqdef \Pi+2(k_1-\varepsilon_1)H^2+2\varepsilon_1(2H^2-K)-4k_1HB$. Then we obtain
\begin{eqnarray*}
&&(T^{\alpha\beta}a_{\beta})_{,\alpha}+(q^{\alpha}\nn)_{,\alpha}\\
&&=-(\Pi{a}^{\alpha\beta}\aab)_{,\alpha}+2\varepsilon_0(S^{\alpha\beta}\aab)_{,\alpha}
+M^{\alpha\mu}b^{\beta}_{\mu,\alpha}\aab+(M^{\alpha\mu}b^{\beta}_{\mu}
b_{\alpha\beta}+M^{\alpha\beta}_{,\alpha\beta})\nn\\
&&=(Pa^{\alpha\beta}\aab)_{,\alpha}+2\varepsilon_0(S^{\alpha\beta}\aab)_{,\alpha}-4k_1Ha^{\alpha\beta}B_{,\alpha}\aab\nonumber\\
&&\qquad+\Big(2k_1(\Delta_{\Gamma}B-2KB)-2(k_1+\varepsilon_1)(\Delta_{\Gamma}H+2H(H^2-K))\Big)\nn.
\end{eqnarray*}
We still use $\Pi$ to denote $P$. Thus, (\ref{eq:velocity-2}) follows easily.

Now we derive the energy law of (\ref{eq:velocity-2}). We infer from (\ref{eq:velocity-2}) that
\begin{eqnarray}
&&\frac{\ud}{\ud{t}}\int\frac{1}{2}\varrho|\vv|^2\ud{S}
=\int\varrho\vv\cdot\frac{\partial{\vv}}{\partial{t}}\ud{S}\nonumber\\
&&=\int\vv\cdot
\Big(-(\Pi{a^{\alpha\beta}\aaa_{\alpha}})_{,\beta}+2\varepsilon_0
(S^{\ab}\aaa_{\alpha})_{,\beta}-4k_1Ha^{\ab}B_{,\beta}\aaa_{\alpha}
\nonumber\\
&&\quad+2k_1\big(\Delta_{\Gamma}{B}+B(b^{\alpha}_{\beta}b^{\beta}_{\alpha}-4H^2)\big)\nn-2\mu_1\big(\Delta_{\Gamma}{H}
+H(b^{\alpha}_{\beta}b^{\beta}_{\alpha}-2H^2)\big)\nn\Big)\ud{S}\nonumber\\
&&=\int\Pi\aaa^{\alpha}\cdot\vv_{,\alpha}-2\varepsilon_0S^{\alpha\beta}\aaa_{\alpha}
\cdot\vv_{,\beta}
-4k_1Ha^{\alpha\beta}B_{,\beta}v_{\alpha}\nonumber\\
&&\quad+2k_1\big(\Delta_{\Gamma}{B}+B(b^{\alpha}_{\beta}b^{\beta}_{\alpha}-4H^2)\big)v^n-2\mu_1\big(\Delta_{\Gamma}{H}
+H(b^{\alpha}_{\beta}b^{\beta}_{\alpha}-2H^2)\big)v^n\ud{S}\nonumber\\
&&=\int-2\varepsilon_0{S}^{\alpha\beta}S_{\alpha\beta}
+2k_1B\Big(2(Hv^{\alpha})_{,\alpha}+(b^{\alpha}_{\beta}b^{\beta}_{\alpha}
-4H^2)v^{n}+\Delta_{\Gamma}v^n\Big)\nonumber\\\label{elawv}
&&\quad-2\mu_1H\Big(\Delta_{\Gamma}{v^n}+(b^{\alpha}_{\beta}b^{\beta}_{\alpha}-2H^2)v^n\Big)\ud{S},
\quad \mu_1=k_1+\veps_1.\label{eq:energy-v}
\end{eqnarray}
Here we used $v^{\alpha}_{,\alpha}-2Hv^n=0$ such that \beno
\aaa^\al\cdot\vv,_\al=v_{,\al}^\al-v^nb^\al_\al=v^{\alpha}_{,\alpha}-2Hv^n=0.
\eeno On the other hand, the Helfrich energy can be simplified as
\beno E_H=\int_{\Gamma}4k_1(H-B)^2+4\varepsilon_1(H^2-K)\ud{S},
\eeno where
$K=\frac{1}{2}(4H^2-b^{\alpha}_{\beta}b^{\beta}_{\alpha})$ is the
Gaussian curvature. As $\int{K}\ud{S}$ is a constant independent
of the time, we have
\begin{eqnarray}\label{eq:energy-Helf}
\frac{\ud}{\ud{t}}E_H&=&\int_{\Gamma}\big(8k_1(H-B)+8\varepsilon_1H\big)
\frac{\partial{H}}{\partial{t}}\ud{S}\nonumber\\\label{elawh}
&=&4\int(\mu_1H-k_1B)(\Delta_{\Gamma}v^n+v^nb^{\alpha}_{\beta}b^{\beta}_{\alpha}+2v^{\alpha}H_{,\alpha}).
\end{eqnarray}
Adding up (\ref{eq:energy-v}) and (\ref{eq:energy-Helf}), and using
$v^{\alpha}_{,\alpha}-2Hv^n=0$ again, we obtain the following energy
dissipation law: \beno \frac{1}{2}\frac{\ud}{\ud
t}\Big(E_H+\int_{\Gamma}\varrho|\vv|^2\ud{S}\Big)=-2\varepsilon_0
\int_{\Gamma}S^{\alpha\beta}S_{\alpha\beta}\ud S. \eeno

\subsection{Some basic estimates in Sobolev spaces}

Let us first recall some product estimates and commutator estimates.

\begin{Lemma}\label{lem:product}
Let $s\ge 0$. Then for any multi-index $\al, \be$, it holds that
\beno
\|\p^\al f\p^\be g\|_{H^s}\le C\big(\|f\|_{L^\infty}\|g\|_{H^{s+|\al|+|\be|}}
+\|g\|_{L^\infty}\|f\|_{H^{s+|\al|+|\be|}}\big).
\eeno
In particular, we have
\beno
\|fg\|_{H^s}\le C\big(\|f\|_{L^\infty}\|g\|_{H^{s}}+\|g\|_{L^\infty}\|f\|_{H^{s}}\big).
\eeno
\end{Lemma}

\begin{Lemma}\label{lem:composition}
Let $s\ge 0$ and $F(\cdot)\in C^\infty(\RR^+)$ with $F(0)=0$. Then
\beno
\|F(f)\|_{H^s}\le C(\|f\|_{L^\infty})\|f\|_{H^s}.
\eeno
\end{Lemma}

\begin{Lemma}\label{lem:commutator}
Let $s>0$. It holds that
\beno
\|\big[\Lam^s, g\big]f\|_{L^2}\le C\big(\|\na g\|_{L^\infty}\|f\|_{H^{s-1}}
+\|g\|_{H^{s}}\|f\|_{L^\infty}\big).
\eeno
Here $\Lam=(-\Delta)^\f 1 2$.
\end{Lemma}

Lemmas \ref{lem:product}-\ref{lem:commutator} are well-known, see \cite{Kato, Triebel} for example.

\begin{Lemma}\label{lem:operator-upper bound}
Let $s\ge 0$ and $k\ge 1$ be an integer. Then it holds that
\beno
\|(a\Delta)^kf\|_{H^s}\le C(\|a\|_{H^2})\big(\|f\|_{H^{s+2k}}+
\|a\|_{H^{s+2k}}\|f\|_{H^2}\big).
\eeno
\end{Lemma}

\no{\bf Proof.}\,We will prove it by induction on $k$. For $k=1$, using Lemma \ref{lem:product}
and the Sobolev inequality, we get
\beno
\|a\Delta f\|_{H^s}&\le& C\big(\|a\|_{L^\infty}\|f\|_{H^{s+2}}+\|a\|_{H^{s+2}}\|f\|_{L^\infty}\big)\\
&\le& C\big(\|a\|_{H^2}\|f\|_{H^{s+2}}+\|a\|_{H^{s+2}}\|f\|_{H^2}\big).
\eeno
Assume Lemma \ref{lem:operator-upper bound} holds  for $k-1$. Then using the induction assumption, we have
\beno
\|(a\Delta)^kf\|_{L^2}&=&\|(a\Delta)^{k-1}(a\Delta f)\|_{L^2}\\
&\le& C(\|a\|_{H^2})\big(\|a\Delta f\|_{H^{s+2k-2}}+\|a\|_{H^{s+2k-2}}\|a\Delta f\|_{H^2}\big).
\eeno
We get by Lemma \ref{lem:product} that
\beno
\|a\Delta f\|_{H^{s+2k-2}}&\le& C\big(\|a\|_{H^2}\|f\|_{H^{s+2k}}+\|a\|_{H^{s+2k}}\|f\|_{H^2}\big),\\
\|a\|_{H^{s+2k-2}}\|a\Delta f\|_{H^2}
&\le& C\|a\|_{H^{s+2k-2}}\big(\|a\|_{H^4}\|f\|_{H^2}+\|a\|_{H^2}\|f\|_{H^4}\big)\\
&\le& C(\|a\|_{H^2})\big(\|f\|_{H^{s+2k}}+\|a\|_{H^{s+2k}}\|f\|_{H^2}\big).
\eeno
Here we used the following interpolation inequality in the last inequality:
\beno
&&\|a\|_{H^4}\le \|a\|_{H^2}^\theta\|a\|_{H^{s+2k}}^{1-\theta},\quad
\|a\|_{H^{s+2k-2}}\le  \|a\|_{H^2}^{1-\theta}\|a\|_{H^{s+2k}}^{\theta},
\eeno
with $\theta=(s+2k-4)/(s+2k-2)$. Thus, we get
\beno
\|(a\Delta)^kf\|_{H^s}\le C(\|a\|_{H^2})\big(\|f\|_{H^{s+2k}}+
\|a\|_{H^{s+2k}}\|f\|_{H^2}\big).
\eeno
The proof is completed. \endproof

\begin{Lemma}\label{lem:operator-lowerbound}
Let $s\ge 0$ and $s_0\in (1,2)$. Assume that $a\ge c_0$ for some positive constant $c_0$. Then we have
\beno
\|(a\Delta)^kf\|_{H^s}\ge c\|f\|_{H^{s+2k}}-C(\|a\|_{H^{s_0+1}})\|a\|_{H^{s+2k}}\|f\|_{H^{s_0}}.
\eeno
\end{Lemma}

\no{\bf Proof.}\,We prove the lemma based on the induction assumption on $k$. For $k=1$, we have
\beno
\|(a\Delta)^kf\|_{H^s}\ge c_0\||\Lam^{s+2}\|_{L^2}-\|\big[\Lam^s,a\big]\Delta f\|_{L^2}.
\eeno
We write
\beno
\big[\Lam^s,a\big]\Delta f=\big[\Lam^s\Delta,a\big]f-2\Lam^s(\na a\cdot\na f)-\Lam^s(\Delta a f),
\eeno
which along with Lemma \ref{lem:product} and Lemma \ref{lem:commutator} implies that
\beno
\|\big[\Lam^s,a\big]\Delta f\|_{L^2}
\le C\big(\|a\|_{H^{s_0+1}}\|f\|_{H^{s+1}}+\|a\|_{H^{s+2}}\|f\|_{H^{s_0}}\big).
\eeno
This yields the case of $k=1$ by an interpolation argument.

Now let us assume that Lemma \ref{lem:operator-lowerbound} holds for $k-1$. Using the induction assumption, we have
\beno
\|\big(a\Delta\big)^kf\|_{L^2}&=&\|(a\Delta)^{k-1}(a\Delta f)\|_{L^2}\\
&\ge& c\|a\Delta f\|_{H^{s+2k-2}}-C(\|a\|_{H^{s_0+1}})\|a\|_{H^{s+2(k-1)}}\|a\Delta f\|_{H^s}.
\eeno
Using the case of $k=1$, we get
\beno
\|a\Delta f\|_{H^{s+2k-2}}\ge c\|f\|_{H^{2k+s}}-C(\|a\|_{H^{s_0+1}})\|a\|_{H^{s+2k}}\|f\|_{H^{s_0}},
\eeno
and by Lemma \ref{lem:product},
\beno
&&\|a\|_{H^{s+2(k-1)}}\|a\Delta f\|_{H^s}\\
&&\le C\|a\|_{H^{s+2(k-1)}}\big(\|a\|_{H^{s_0+2}}\|f\|_{H^{s_0}}+\|a\|_{H^{s_0}}\|f\|_{H^{s_0+2}}\big)\\
&&\le \veps\|f\|_{H^{s+2k}}+C(\|a\|_{H^{s_0+1}})\|a\|_{H^{s+2k}}\|f\|_{H^{s_0}}.
\eeno
Here we used the following interpolation inequality in the last inequality:
\beno
&&\|a\|_{H^{s_0+2}}\le \|a\|_{H^{s_0}}^\theta\|a\|_{H^{s+2k}}^{1-\theta},\quad
\|a\|_{H^{s+2k-2}}\le  \|a\|_{H^{s_0}}^{1-\theta}\|a\|_{H^{s+2k}}^{\theta},
\eeno
with $\theta=2/(s+2k-s_0)$. Taking $\veps$ to be small enough, we obtain
\beno
\|(a\Delta)^kf\|_{H^s}\ge c\|f\|_{H^{s+2k}}-C(\|a\|_{H^{s_0+1}})\|a\|_{H^{s+2k}}\|f\|_{H^{s_0}}.
\eeno
The proof is completed.\endproof

\begin{Lemma}\label{lem:operator-commuatator-time}
Let $s\ge 0$ and $k\ge 1$ be an integer. Then we have
\beno
&&\|\big[\p_t, (a\Delta)^k\big]f\|_{H^s}\le
C(\|(a,\p_ta)\|_{H^2})\big(\|f\|_{H^{s+2k}}+ \|(a,\p_t
a)\|_{H^{s+2k}}\|f\|_{H^2}\big),\\
&&\|\big[\na, (a\Delta)^k\big]f\|_{H^s}\le C(\|a\|_{H^3})\big(\|f\|_{H^{s+2k}}+
\|a\|_{H^{s+2k+1}}\|f\|_{H^2}\big).
\eeno
\end{Lemma}

\no{\bf Proof.}\,We write
\beno
\big[\p_t, (a\Delta)^k\big]f&=&\sum_{\ell=1}^{k-1}(a\Delta)^\ell\big[\p_t, a\Delta\big](a\Delta)^{k-\ell-1}f\\
&=&\sum_{\ell=1}^{k-1}(a\Delta)^\ell(\p_ta\Delta)(a\Delta)^{k-\ell-1}f.
\eeno
Then the first inequality can be seen easily from the proof of Lemma \ref{lem:operator-upper bound}.
The proof of the second inequality is similar.\endproof

\begin{Lemma}\label{lem:operator-commuatator-function}
Let $s\ge 0$ and $k\ge 1$ be an integer. Then we have
\beno
\|\big[(a\Delta)^k, g\big]f\|_{H^s}\le C(\|a\|_{H^2},\|g\|_{H^3})\big(\|f\|_{H^{s+2k-1}}+
\|(a,g)\|_{H^{s+2k}}\|f\|_{H^2}\big).
\eeno
\end{Lemma}

\no{\bf Proof.}\,As in Lemma \ref{lem:operator-lowerbound}, this lemma can be
proved by the induction argument, however, we omit the details here.\endproof

\subsection{Elliptic estimates}

We consider the following elliptic system:
\begin{eqnarray}\label{eq:elliptic system}
\left\{
\begin{array}{l}
\Big(\frac{W_1}{\se}\Big)_{\alpha}-\Big(\frac{W_2}{\se}\Big)_{\beta}=f_1,\\
\Big(\frac{W_1}{\se}\Big)_{\beta}+\Big(\frac{W_2}{\se}\Big)_{\alpha}=f_2.
\end{array}\right.
\end{eqnarray}
We write
\beno
\f {W_1} {\sqrt{E}}=\p_\al\phi+\p_\be\psi,\quad \f {W_2} {\sqrt{E}}=-\p_\be\phi+\p_\al\psi.
\eeno
Then (\ref{eq:elliptic system}) is reduced to solve the following Poisson equations:
\beno
\Delta\phi=f_1,\quad \Delta\psi=f_2.
\eeno
Thus, we have
\begin{Lemma}\label{lem:elliptic system}
Let $s\ge 1$. If $E\in H^{s}(\RT^2), f_1,f_2\in H^{s-1}(\RT^2)$, then the system (\ref{eq:elliptic system}) has
a solution $(W_1,W_2)$ satisfying
\beno
\|W_1\|_{H^{s-1}}+\|W_2\|_{H^{s-1}}\le C(\|E\|_{H^s})\big(\|f_1\|_{H^{s-1}}+\|f_2\|_{H^{s-1}}\big).
\eeno
\end{Lemma}

Next we consider the elliptic equation:
\ben\label{eq:elliptic equation}
-\Delta U+a U=f.
\een
\begin{Lemma}\label{lem:elliptic equation}
Let $s\ge 0$. Assume that $f\in H^s(\RT^2)$, and $a\in H^{s}\cap L^\infty(\RT^2)$ with
\ben\label{eq:curvature condition}
a(x)\ge 0,\quad \int_{\RT^2}a(x)dx\ge a_0>0.
\een
Then there exists a unique solution $U\in H^{s+2}(\RT^2)$ to (\ref{eq:elliptic equation}) satisfying
\beno
\|U\|_{H^{s+2}}\le C\|f\|_{H^s}.
\eeno
Here $C$ is a constant depending only on $a_0$ and $\|a\|_{ H^{s}\cap L^\infty}.$
\end{Lemma}

\no{\bf Proof.\,}The proof of existence part is standard. Here we only prove the  estimate.
Taking the $L^2$ inner estimate gives
\beno
\|\na U\|_{L^2}^2+\int_{\RT^2}a|U|^2dx=\int_{\RT^2}fUdx.
\eeno
Let $\bar U=\f 1 {4\pi^2}\int_{\RT^2}Udx$. Therefore, we have
\beno
\|U\|_{L^2}^2\le 2\|U-\bar U\|_{L^2}^2+2\|\bar U\|_{L^2}^2\le 2\|\na U\|_{L^2}+8\pi^2|\bar U|.
\eeno
On the other hand, we have by (\ref{eq:curvature condition}) that
\beno
a_0|\bar U|^2\le \int_{\RT^2} a|\bar U|^2dx&\le& 2\int_{\RT^2} a|U|^2dx
+\int_{\RT^2}a|U-\bar U|^2dx\\
&\le& 2\int_{\RT^2} a|U|^2dx+\|a\|_{L^\infty}\|\na U\|^2_{L^2}.
\eeno
This yields that
\ben\label{eq:L2-est}
\|U\|_{L^2}\le C\big(\|\na U\|_{L^2}^2+\int_{\RT^2}a|U|^2dx\big).
\een
Using the elliptic estimate in $H^s$, we obtain
\beno
\|U\|_{H^{s+2}}\le C\big(\|aU\|_{H^s}+\|f\|_{H^s}\big)
\le C\big(\|a\|_{L^\infty}\|U\|_{H^s}+\|a\|_{H^s}\|U\|_{L^\infty}+\|f\|_{H^s}\big),
\eeno
from which and (\ref{eq:L2-est}), the desired estimate follows from an interpolation argument.
\endproof

\bigskip

\noindent {\bf Acknowledgments.} The authors are grateful to the referees and the editor for their
invaluable comments and suggestions which have helped us improve the paper significantly.
The authors are grateful too to
Jianzhen Qian, Dan Hu, and Peng Song for sharing their sights in many discussions.
Wei Wang and Pingwen Zhang are supported by
the NSF of China under Grant 11011130029.
Zhifei Zhang is supported by the NSF of China under Grants 10990013 and 11071007.


\end{document}